{
\documentclass[envcountsame,envcountsect]{svmult}

\usepackage{mathptmx}       
\usepackage{helvet}         
\usepackage{courier}        
\usepackage{type1cm}        
\usepackage{makeidx}         
\usepackage{graphicx}        
\usepackage{multicol}        
\usepackage[bottom]{footmisc}

\usepackage{listings}

\usepackage{booktabs}

\usepackage{hyperref}

\makeindex             

\usepackage{amstext,amsfonts,amssymb,amscd,amsbsy,amsmath}

\smartqed

\newcommand{\A}{\mathbb{A}}
\newcommand{\PP}{\mathbb{P}}
\newcommand{\ZZ}{\mathbb{Z}}

\newcommand{\vh}{\vec{h}}

\DeclareMathOperator{\Gr}{Gr}
\DeclareMathOperator{\Hom}{Hom}
\newcommand{\cH}{\mathcal{H}}

\DeclareMathOperator{\lead}{lead}
\DeclareMathOperator{\Span}{span}
\DeclareMathOperator{\Proj}{Proj}

\DeclareMathOperator{\init}{in}

\newcommand{\sw}[1]{\emph{#1}} 

    \newcommand{\kk}{\Bbbk}
    \DeclareMathOperator{\Spec}{Spec}
\newcommand{\hook}{\ensuremath{\mathbin{\text{\raisebox{.4ex}{%
  \vrule height .5pt width 1ex depth 0pt%
  \vrule height 0.8ex width .5pt depth 0pt%
}}\mathchoice{}{}{\mkern3mu}{\mkern3mu}}}}%

\begin{document}
\title*{The Hilbert scheme of $11$ points in $\A^3$ is irreducible}


\author{Theodosios Douvropoulos, Joachim Jelisiejew,
    Bernt Ivar Utst{\o}l N{\o}dland, and Zach Teitler}
\authorrunning{T. Douvropoulos, J. Jelisiejew, B.I. U.
    N{\o}dland, Z. Teitler}
\institute{Theodosios Douvropoulos \at School of Mathematics, University of Minnesota, Minneapolis, \email{douvr001@umn.edu}
\and Joachim Jelisiejew \at Faculty of Mathematics, Informatics and Mechanics,
University of Warsaw, Poland \email{jjelisiejew@mimuw.edu.pl}
\and Bernt Ivar Utst{\o}l N{\o}dland \at Department of Mathematics, University of Oslo, Norway. \email{berntin@math.uio.no}
\and Zach Teitler \at Boise State University, Department of Mathematics, 1910 University Drive, Boise, Idaho 83725-1555, USA. \email{zteitler@boisestate.edu}}

\date{\today}

%
%
\maketitle


\abstract{
    We prove that the Hilbert scheme of $11$ points on a smooth threefold is
    irreducible. In the course of the proof, we present several known and new
    techniques for producing curves on the Hilbert scheme.
}


\section{Introduction}\label{section: Introduction}

    Let $X$ be a smooth connected quasi-projective variety.
    The Hilbert scheme of $d$ points in $X$ is the scheme
    parametrizing finite subschemes of $X$ of degree $d$.
    There are ample introductory readings on the Hilbert scheme of points
    available, including~\cite{fantechi, fogarty,
    Gottsche_Hilbert_schemes_and_Betti_numbers, hartshorne_deformation_theory,
iakanev, Miller_Sturmfels, nakajima_lectures_on_Hilbert_schemes}.

    The Hilbert scheme of points is quasi-projective (projective iff $X$ is) and
    connected~\cite{fantechi, fogarty, FGA}.
    Moreover, Fogarty~\cite{fogarty} proved that for $\dim X \leq 2$ it is
    smooth of dimension $d\cdot(\dim X)$.
    For higher-dimensional $X$, much less is known.
    The questions of irreducibility of
    the Hilbert scheme of points is especially interesting, because it ensures
    that all finite schemes are limits of reduced ones;
    see~\cite{bubu2010} for an application.
    This question is local and only depends on the dimension of $X$:
    the answer for $n$-dimensional $X$ will be the same as for $\A^n$,
    see~\cite[p.~4]{artin_deform_of_sings} or \cite[Lemma~2.2]{cn09}.
    We denote the Hilbert scheme of $d$ points in $\A^n$ by
    $H_{n}^d$. Our motivating question is the following:
    \begin{center}
        \emph{For which pairs $(n, d)$ is the Hilbert scheme $H_{n}^d$ irreducible?}
    \end{center}
    By Fogarty's results, all $H_2^d$ are irreducible.
    Mazzola~\cite{Mazzola} proved irreducibility of
    $H_{n}^d$ for all $n$ and $d\leq 7$.
    Iarrobino~\cite{iarrobino_reducibility, iarrobino_compressed_artin} showed that for every $n\geq 3$
    and ${d\geq 78}$ the scheme $H_{n}^d$ is reducible.
    Emsalem and
    Iarrobino proved that
    $H_n^d$ is reducible for $d\geq 8$ and $n\geq 4$, see \cite[Section 2.2,
    p.~158]{emsalem_iarrobino_small_tangent_space} and also~\cite{CEVV}.
    Borges dos Santos, Henni, and Jardim~\cite{dosSantos}  showed that
    $H_{3}^9$ and $H_{3}^{10}$ are irreducible by comparing
    them with appropriate spaces of commuting matrices and using the results of {\v{S}}ivic~\cite[Theorems 26,~32]{Sivic__Varieties_of_commuting_matrices}.
    Thus, the reducibility of $H_{n}^d$ was unknown only for the values
    $n = 3$ and $11\leq d\leq 77$. Here we improve the lower bound.

    \begin{theorem}\label{ref:main:theorem}
        The Hilbert scheme of $11$ points in a smooth irreducible threefold is
        irreducible of dimension $33$.
    \end{theorem}

    We prove Theorem~\ref{ref:main:theorem} in Section~\ref{section: Cases}.
    We review background information in Section~\ref{section: Prerequisites}.
    In Section~\ref{section: Preparation} we give an overview of strategy,
    gather general results that will be used in the proof of the above theorem,
    and demonstrate how to use \sw{Macaulay2} \cite{M2} for some computations.

    In Section~\ref{sec:upto95} we discuss a special class of subschemes,
    which appeared in the earliest example of reducible $H_{3}^d$, due to
    Iarrobino~\cite{iarrobino_reducibility}. Namely, let $\mathfrak{m}$ be the ideal
    of the origin of $\A^3$. Fix
    $d$ and  consider the ideals $\mathfrak{m}^s \subset I \subset
    \mathfrak{m}^{s+1}$ such that $V(I)$ has degree $d$; then $s$ is
    uniquely determined. Call such ideals \emph{very compressed} and denote by
    $\cH^{\max,d}$ their family.
    Let $R_3^d$ denote the closure in $H_3^d$ of the
    open set of smooth subschemes. The component $R_3^d$ is called
    the \emph{smoothable component}. It has dimension $3d$.
    The result of~\cite{iarrobino_reducibility} is that for $d\geq 96$ we have
    $\dim \cH^{\max, d} \geq 3d$ and, thus, a general very compressed
    ideal does not lie in the smoothable component.
    We show that for $d\leq 95$ the family $\cH^{\max, d}$ is in fact
    contained in $R^{d}_3$.
    \begin{proposition}\label{ref:upto95:prop}
        The family $\cH^{d, \max}$ of very compressed ideals is
        contained in the smoothable component if and only if $d\leq 95$.
    \end{proposition}
    The key points of the proof are the use of smoothings by degenerating to
    initial ideals and a \sw{Macaulay2} calculation, see
    Section~\ref{sec:upto95}.

We now explain our approach to the proof of Theorem~\ref{ref:main:theorem}.
We build upon the strategy of~\cite{CEVV}.
As explained there, questions about smoothability of a specified ideal
$I$ are easily reduced to the case where $I$ is local and has full embedding dimension~$3$.
There are fifteen possible Hilbert functions of $I$,
see Table~\ref{table: case proof sections}.
For each Hilbert function $\vh$,
the scheme $H_3^{\vh}$ parameterizes local ideals with fixed Hilbert function $\vh$
and the \emph{standard graded Hilbert scheme} $\cH_3^{\vh}$ parameterizes homogeneous ideals with fixed Hilbert function $\vh$.
We apply three different strategies to show that for each Hilbert function
$\vh$ in our list, we have $H_3^{\vh}
\subset R_3^{11}$.

\begin{table}
\begin{tabular}{ll@{\hspace{1.5em}}l @{\hspace{2.5em}} ll@{\hspace{1.5em}}l @{\hspace{2.5em}} ll@{\hspace{1.5em}}l}
\toprule
1. & $(1,3,1,1,1,1,1,1,1)$ & \S\ref{sec:tailsofones} &
6. & $(1,3,5,1,1)$ & \S\ref{sec:tailsofones} &
11. & $(1,3,2,2,2,1)$ & \S\ref{sec:132221} \\
2. & $(1,3,2,1,1,1,1,1)$ & \S\ref{sec:tailsofones} &
7. & $(1,3,3,4)$ & \S\ref{sec:1def} &
12. & $(1,3,3,2,2)$ & \S\ref{sec:13322} \\
3. & $(1,3,2,2,1,1,1)$ & \S\ref{sec:tailsofones} &
8. & $(1,3,4,3)$ & \S\ref{sec:1def} &
13. & $(1,3,3,3,1)$ & \S\ref{sec:13331} \\
4. & $(1,3,3,1,1,1,1)$ & \S\ref{sec:tailsofones} &
9. & $(1,3,5,2)$ & \S\ref{sec:1def} &
14. & $(1,3,3,2,1,1)$ & \S\ref{sec:133211} \\
5. & $(1,3,4,1,1,1)$ & \S\ref{sec:tailsofones} &
10. & $(1,3,4,2,1)$ & \S\ref{sec:13421} &
15. & $(1,3,6,1)$ & \S\ref{sec:1361} \\
\bottomrule
\end{tabular}
\caption{Hilbert functions $\vh$ in $H_3^{11}$ and the corresponding sections.}\label{table: case proof sections}
\end{table}

First, for some cases the knowledge about the Hilbert function of an ideal
$I$ is enough to produce a deformation (via \emph{ray families} introduced
in~\cite{cjn13}) whose special fiber is $I$ and general fiber is reducible.
By Lemma~\ref{lemma: cleavable implies smoothable}, such an $I$ is smoothable,
see Section~\ref{sec:tailsofones}.

Second, most of the schemes $H^{\vh}_3$ contain smooth points
of the Hilbert scheme which lie in the smoothable component $R^{11}_3$.
Such points are called
\emph{smooth and smoothable points}; examples include points corresponding to
Gorenstein algebras, see~\cite[Corollary~2.6]{cn09}.
\begin{lemma}
\label{lemma: smooth and smoothable is enough}
If $Z \subseteq H_3^{11}$ is an irreducible set
that contains a smooth and smoothable point,
then we have $Z \subseteq R_3^{11}$.
\end{lemma}
\begin{proof}
    The locus of smooth and smoothable points is open and contained in
    $R^{11}_3$,
    so the intersection $Z \cap R^{11}_{3}$ contains an open subset of $Z$.
    Then, the subset $Z\cap
    R^{11}_3 \subset Z$ is dense and closed, so it is equal to $Z$.\qed
\end{proof}
To apply the above lemma, we write $H_3^{\vh}$ as a union of irreducible sets
$Z$ and show that each $Z$ contains a smooth and smoothable point.
To find the sets $Z$ we may take advantage of the
morphism $\pi_{\,\vh}:H^{\vh}_3\to\cH^{\vh}_3$
taking an ideal $I$ to its initial ideal, see~\cite{CEVV}.
We employ the following 3-step strategy:
\begin{enumerate}
\item \label{coarsecomponents}
Decompose $\cH^{\vh}_3$ into irreducible strata.
\item \label{finecomponents}
    Using the morphism $\pi_{\vh}:H^{\vh}_3\to\cH^{\vh}_3$, decompose
    $H^{\vh}_3$ into irreducible strata.
\item \label{contained?}
For each stratum of $H^{\vh}_3$, find a smooth point of the Hilbert scheme
which lies in the smoothable component and conclude that the whole stratum
lies there.
\end{enumerate}

In steps \ref{coarsecomponents} and \ref{finecomponents},
we use Macaulay's inverse systems, see Section~\ref{section: Prerequisites}.
In the simplest cases, we find that there is a bijection between
irreducible strata of $\cH^{\vh}_3$ and $H^{\vh}_3$, but this is not always
true, see for example Section~\ref{sec:13322}.

For step~\ref{contained?} we introduce \emph{cleavable}
ideals. An ideal is said to be \emph{cleavable} (or \emph{limit-reducible})
if it can be deformed to an ideal whose support consists of at least two points.
\begin{lemma}\label{lemma: cleavable implies smoothable}
    A cleavable ideal $I\in H_{3}^{11}$ is smoothable.
\end{lemma}
\begin{proof}
Let $I_t$ be a one-parameter flat family of ideals with $I_0 = I$ and for $t \neq 0$, $I_t$ supported at more than one point.
Each irreducible component of $I_t$ has length strictly less than $11$, so
it is smoothable. Hence, the ideal $I$ is also smoothable.
\qed
\end{proof}
To show that an ideal
$I$ is cleavable, we construct a family over $\Spec \kk[t]$ whose general fiber is reducible and
check that it is flat, see Section~\ref{section: Macaulay2 code examples}.

Third, there is a case where both previous methods do not apply. This is the case $\vh
= (1, 3, 6, 1)$, see Proposition~\ref{proposition: 1361}. The stratum $H^{\vh}_{3}$ does
not seem to contain smooth points. However, the stratum is irreducible and we
can describe what general points look like. We build a deformation showing
that such general points are smoothable, hence, by irreducibility, the entire
stratum has to be smoothable.

We work over an algebraically closed field $\kk$ of characteristic zero.


\section{Prerequisites}\label{section: Prerequisites}

\runinhead{Hilbert schemes and smoothability.}
The \emph{Hilbert scheme} $H^d_n$ parameterizes subschemes
of $\A^n$ of dimension zero and degree $d$.
More formally, $H^d_n$ represents the functor which assigns to
each $\kk$-scheme $X$ the set of subschemes of $\A^n \times X$
which are flat over $X$ and for which all fibers are finite of degree $d$,
see~\cite[Chapter~1]{hartshorne_deformation_theory}.
Equivalently, letting $T=\kk[\alpha_1,\alpha_2, \dots,\alpha_n]$,
the scheme $H^d_n$ parameterizes ideals $I$ for which $T/I$
is a vector space of dimension $d$.
In other words, $H^d_n$ also represents the functor which assigns to each $\kk$-algebra $A$
the set of ideals $I$ in $T \otimes A$ such that the quotients $T\otimes A/I$ are locally free $A$-modules
of rank $d$.
\par
The Zariski tangent space to $H^d_n$ at the point representing $I$ is the
$T$-module
$\Hom(I,T/I)$, see~\cite[Theorem 1.1]{hartshorne_deformation_theory}.
Using \sw{Macaulay2}~\cite{M2}, we can compute the dimension of this tangent
space.
We stress that a point is smooth if and only if the point lies on only one irreducible component of the scheme
and the dimension of the tangent space at that point
equals the dimension of the component of the scheme containing the point.
The dimension of the tangent space increases at singular points.

On $H_{n}^d$, there is a distinguished component corresponding to smooth
schemes. Indeed, a slightly perturbed tuple of $d$ closed points in
$\A^n$ is just another such tuple. Thus, the set of tuples of points is open in
the Hilbert scheme and their closure is a component. It is called the
\emph{smoothable component} of $H_{n}^d$ and denoted by
$R_{n}^d$.
Clearly, $R_{n}^d$ is generically smooth of dimension
$nd$. Since $H_2^d$ is smooth, we have $R_{2}^d = H_{2}^d$.

A point of $R_{n}^d$ is said to be \emph{smoothable}.
Thus, an ideal $I$ is smoothable if and only if it can be deformed to an ideal of $d$ distinct points.
This means that
one can build a one-parameter flat family of schemes
over a discrete valuation ring
for which the general member consists of $d$ distinct points and the special fiber
is $T/I$, see~\cite{jabu_jelisiejew_smoothability, CEVV} for details.
In particular, a disjoint union of smoothable schemes is smoothable and a limit
of smoothable schemes is smoothable.

\runinhead{Hilbert functions.}
In analyzing the Hilbert scheme $H^d_n$, it is useful to use work with an invariant that refines the degree $d$.
There are two closely-related notions of \emph{Hilbert function}:
\begin{itemize}
\item
For a graded $T$-module $M$,
its Hilbert function is defined by $\vh(i) = \dim(M_i)$.
In particular, given a homogeneous ideal $I \subset T$,
we consider the Hilbert function of the quotient ring $T/I$.
\item
For a filtered $T$-module $M$
with descending filtration $M = M_0 \supseteq M_1 \supseteq M_2 \supseteq \dotsb$,
the Hilbert function $\vh$ is defined by $\vh(i) = \dim(M_i / M_{i+1})$.
In particular, if the scheme associated to an ideal $I \subset T$ is supported at a point,
then $T/I$ is a local ring $(A,\mathfrak{m})$, and the Hilbert function $\vh$
with respect to the filtration by powers of $\mathfrak{m}$ is defined to be
\(
  \vh(i)=\dim(\mathfrak{m}^i/\mathfrak{m}^{i+1}).
\)
If $I$ is homogeneous and $T/I$ is local, the two notions coincide.
\end{itemize}
We write $\vh$ as a vector $(\vh(0),\vh(1),\dots)$,
trimming it after the last positive entry.

Let $A = T/I$ where $T = \kk[\alpha_1,\alpha_2,\dotsc,\alpha_n]$ is a
polynomial ring with its standard grading and $I$ is a homogeneous ideal. Assume that $I$ contains no
linear forms.  We call such an algebra \emph{standard graded}.

Macaulay's bound is an upper bound for the growth of Hilbert functions of standard graded algebras,
defined as follows.
First, for positive integers $h$ and $d$, there exist uniquely determined integers
$\delta \geq 1$ and $k_d > k_{d-1} > \dotsb > k_\delta \geq \delta$ such that
\begin{equation*}
  h = \binom{k_d}{d} + \binom{k_{d-1}}{d-1} + \dotsb + \binom{k_\delta}{\delta} .
\end{equation*}
This expression is called the \emph{$d$-binomial expansion} of $h$ and denoted $h_{(d)}$.
The $d$\nobreakdash-binomial expansion of $h$ can be found greedily: let $k_d$ be the greatest integer such that
$\binom{k_d}{d} \leq h$, then find the $(d-1)$-binomial expansion of $h-\binom{k_d}{d}$.
Now $h^{\langle d \rangle}$ is defined as follows.
If $h_{(d)} = \binom{k_d}{d} + \binom{k_{d-1}}{d-1} + \dotsb +
\binom{k_\delta}{\delta}$ then we define
\begin{equation*}
h^{\langle d \rangle} := \binom{k_d+1}{d+1} + \dotsb + \binom{k_\delta+1}{\delta+1}.
\end{equation*}

\begin{example}
    We have $5_{(2)} = \binom{3}{2} + \binom{2}{1}$, so $5^{\langle 2 \rangle}
    = \binom{4}{3} + \binom{3}{2} = 7$.  Similarly, we have $4_{(2)} =
    \binom{3}{2} + \binom{1}{1}$, so $4^{\langle 2 \rangle} = \binom{4}{3} +
    \binom{2}{2} = 5$.
\end{example}

\begin{example}
If $h \leq d$ then we have $h_{(d)} = \binom{d}{d} + \binom{d-1}{d-1} + \dotsb + \binom{d-h+1}{d-h+1}$
and $h^{\langle d \rangle} = h$.
\end{example}

\begin{theorem}[Macaulay's bound, {\cite{macaulay_enumeration} or \cite[Theorem 4.2.10]{BrunsHerzog}}]
\label{Macaulay's bound}
Let $A$ be a standard graded $\kk$-algebra with Hilbert function $\vh$.
For every non-negative integer $d$, we have $\vh(d+1) \leq \vh(d)^{\langle d \rangle}$.
\end{theorem}

\begin{corollary}\label{corollary: hilbert function nonincreasing}
Let $A$ be a standard graded $\kk$-algebra with Hilbert function $\vh$.
If $d \geq 0$ is such that $\vh(d) \leq d$, then we have $\vh(d) \geq \vh(d+1) \geq \vh(d+2) \geq \dotsb$.
\end{corollary}

Once the Macaulay bound is attained then it will also be attained for all
higher degrees provided that no new generators of the ideal appear:
\begin{theorem}[Gotzmann's Persistence Theorem, {\cite{gotzmann_persistence_theorem} or \cite[Theorem~4.3.3]{BrunsHerzog}}]
\label{Gotzmann's Persistence}
Let $A = T/I$ be a standard graded algebra with Hilbert function $\vh$.
If $d \geq 0$ is an integer such that $\vh(d+1) = \vh(d)^{\langle d \rangle}$ and $I$ is generated in degrees $\leq d$,
then we have $\vh(k+1) = \vh(k)^{\langle k \rangle}$ for all $k \geq d$.
\end{theorem}

\runinhead{Apolarity and inverse systems.}
A key tool in the analysis of finite schemes is the technique of \emph{Macaulay's inverse
systems}, also known as apolarity.
General references include \cite{emsalem,
geramita_inverse_systems_of_fat_points}, \cite[Section 1.3, Chapter
5]{iakanev}, \cite{ranestad_schreyer_VSP}.

Let $S = \kk[x_1,x_2,\dotsc,x_n]$
and $T=\kk[\alpha_1,\alpha_2,\dotsc,\alpha_n]$ be polynomial rings with the standard grading.
When $n \leq 3$, we instead use variables $x,y,z$ and $\alpha, \beta, \gamma$.
We write $S_{\leq d}$ for $\bigoplus_{k=0}^d S_k$, and similarly $T_{\leq d}$.
The polynomial ring $T$ acts on $S$ by letting $\alpha_i$ act as partial differentiation by $x_i$.
This is called the \emph{apolarity} action.
We denote this action by $\hook$, so that $\alpha_i \hook F = \tfrac{\partial F}{\partial x_i}$ for $F \in S$.
This gives bilinear maps $T_d \times S_e \to S_{e-d}$ for all $d,e$.
In particular, for each $d$ the pairing $T_d \times S_d \to S_0 = \kk$ is a perfect pairing.

\begin{definition}
For any subset $J \subset S$ the \emph{apolar ideal}, or \emph{annihilating
ideal} $J^\perp \subset T$
is the ideal of elements $\Theta \in T$ such that $\Theta \hook F = 0$ for all $F \in J$.
For $F \in S$ we write $F^\perp$ for $(\{F\})^\perp$.
\end{definition}

When $J$ is spanned by homogeneous elements, the apolar ideal is
homogeneous.  When $J$ consists of a single element $F$, then the ideal
$F^{\perp}$ is \emph{Gorenstein}, see~\cite[Section 21.2]{eisenbud}.

\begin{example}
    If $F = x_1^{a_1} x_2^{a_2}\dotsm x_n^{a_n}$,
then we claim that $F^\perp = (\alpha_1^{a_1+1}, \dotsc, \alpha_n^{a_n+1})$.
Indeed, it is easy to see that each $\alpha_i^{a_i+1} \in F^\perp$.
Conversely, if $\Theta \in T$ has a term $\alpha_1^{b_1}\alpha_2^{b_2} \dotsm \alpha_n^{b_n}$ with each $b_i \leq a_i$,
then the apolar pairing of this term with $F$ is a monomial that determines the $b_i$,
meaning that it cannot be cancelled by the other terms of $\Theta$.
Hence, if $\Theta \in F^\perp$, then each term of $\Theta$ must lie in the indicated ideal.
\end{example}

The linear map $T \to S$ given by $\Theta \mapsto \Theta \hook F$ provides a simple approach to computing $F^\perp$.
The apolar ideal $F^\perp$ is  the kernel of this map.
We can compute $J^\perp$ by intersecting the ideals $F^\perp$ for $F$ in $J$.
If $J$ is a $\kk$-vector space, then it is sufficient to consider $F$ in a basis for $J$.

\begin{example}
For $F = x^3 + yz$, we have $F^\perp = (\alpha^3 - 6 \beta \gamma, \alpha\beta, \alpha\gamma, \beta^2, \gamma^2)$.
\end{example}

\begin{example}
For $F = x^2 y + y^2 z$, we have
$F^\perp = (\gamma^2, \alpha\gamma, \alpha^2 - \beta\gamma, \beta^3, \alpha\beta^2)$.
\end{example}

\begin{definition}
A \emph{Macaulay inverse system}, or simply \emph{inverse system}, is a $T$-sub\-mod\-ule of $S$.
That is, an inverse system is a $\kk$-vector subspace $J \subseteq S$ which is
closed under differentiation: if $F \in J$, then all of the
derivatives ${\alpha_1 \hook F}, \dotsc, {\alpha_n \hook F}$ lie in $J$.
\end{definition}

The \emph{inverse system generated by} a subset $f_1, \ldots , f_s$ of $S$ is
\(
  \langle f_1,f_2,\dotsc,f_s \rangle = Tf_1 + Tf_2 + \dotsb + Tf_s,
\)
 that is, the vector space spanned by the $f_i$ together with all higher partial derivatives.
 Clearly, we have $\langle f_1,\dotsc,f_s \rangle^\perp = \bigcap_{i=1}^s
 \langle f_i\rangle^\perp = \bigcap_{i=1}^s  f_i^\perp$.
An inverse system is \emph{homogeneous} if it is generated by homogeneous
elements.

\begin{remark}
The mapping $J \mapsto J^\perp$ sends finite-dimensional inverse systems
to local ideals supported at the origin, that is, $\mathfrak{m}$-primary ideals
where $\mathfrak{m}$ is the ideal of the origin.
The mapping is one-to-one, since $J$ may be computed from $J^\perp$ similarly to the discussion above.
In fact it is a bijection, as shown by Macaulay \cite{Macaulay_inverse_systems},
or see for example \cite[Corollaire 2]{emsalem}.
When $I$ is a local ideal, we will write $I^\perp$ for its inverse system.

Recall that $H^{\vh}_n$ and $\cH^{\vh}_n$ consist of all homogeneous and local ideals, respectively,
with Hilbert function $\vh$.
On the other hand $H^d_n$, consists of  all zero-dimensional schemes of length $d$ in $\A^n$,
not only local ones or ones supported at the origin.
\end{remark}

\begin{proposition}
[{\cite[Remark after Proposition 2.5]{geramita_inverse_systems_of_fat_points}}]
\label{proposition: homogeneous inverse system isomorphic to apolar algebra}
If $J$ is a homogeneous inverse system then, 
$J$ is isomorphic as a graded $\kk$-vector space to $T/J^\perp$.
\end{proposition}

\begin{proposition}[{\cite[Proposition 2(a)]{emsalem}}]
\label{proposition: inverse system same dimension as apolar algebra}
For a finite-dimensional inverse system $J$, we have
$\dim_{\kk} J = \dim_{\kk} T/J^\perp$.
\end{proposition}

\begin{proof}
Let $d$ be large enough so that $J \subseteq S_{\leq d}$.
It follows that the map $T_{\leq d} \to T/J^\perp$ is surjective.
Hence, both of the dimensions are equal to the codimension of $J^\perp \cap T_{\leq d}$ in $T_{\leq d}$.
\qed
\end{proof}

\begin{remark}

For an inverse system $J$, for each integer $k$,
$J_{\leq k}$ denotes the vector space of polynomials of degree at most $k$ in $J$.
These form an increasing filtration, $J_{\leq 0} \subseteq J_{\leq 1} \subseteq \dotsb$.
The inverse system $J$ is a filtered $T$-module, so its Hilbert function $\vh$ is given by
$\vh(k) = \dim J_{\leq k} - \dim J_{\leq k-1}$ for each $k$ and $\sum \vh(i) = \dim_{\kk} J$.
If $J$ is homogeneous, then $\vh(k) = \dim J_k$.
\end{remark}

\begin{proposition}[{\cite[Lemma~2.12]{iakanev}}]\label{ref:symmetry}
Let $f \in S$ be a homogeneous form of degree $d$. If $\vh$ is the Hilbert function of the inverse system $\langle f \rangle$,
then $\vh = (\vh(0),\dotsc,\vh(d))$ is symmetric: $\vh(i) = \vh(d-i)$ for all $i$.
\end{proposition}

\begin{proposition}[\cite{MR2279854}]\label{proposition: *k1 *11 *221}
Suppose that $f \in S$ is a homogeneous form of degree $d$.
Let $\vh$ be the Hilbert function of $\langle f \rangle$.
If $\vh(d-1) = k$, that is $\vh = (\dotsc,k,1)$,
then there are independent linear functions $\ell_1, \ell_2, \dotsc,\ell_k \in S_1$
and a homogeneous form $g$ such that $f = g(\ell_1,\ell_2,\dotsc,\ell_k)$.
Equivalently, there is a linear change of coordinates so that $f$ depends only on the variables $x_1,\dotsc,x_k$,
not on $x_{k+1},\dotsc,x_n$.
\end{proposition}
\begin{remark} \label{remark: *11 *221}
Using the above proposition, one can show that if $\langle f \rangle$ has Hilbert function $(\dotsc,1,1)$,
then $f = \ell^d$ for some linear function $\ell$ and $\langle f \rangle$ has Hilbert function $(1,1,\dotsc,1,1)$.
If $\vh(d-2) = \vh(d-1) = 2$,
then either $f = \ell^d + m^d$ or $f = \ell^{d-1} m$ for some independent linear functions $\ell, m \in S_1$,
and either way $\langle f \rangle$ has Hilbert function $(1,2,2,\dotsc,2,2,1)$. For proof see for example \cite[Theorem 1.44]{iakanev}:
in their notation, $s=2$, and $f^\perp$ has a quadratic generator,
which up to a change of coordinates is either $\alpha\beta$ or $\beta^2$.
\end{remark}

Dealing with nonhomogeneous inverse systems is much harder than working with
homogeneous ones. Fortunately, each inverse system $J$ has an associated
homogeneous inverse system $\lead(J)$.

\begin{definition}
The \emph{leading form} of a polynomial is its highest degree homogeneous part.
This may not be a monomial.
For an inverse system $J \subset S$, the inverse system of leading forms of $J$, denoted $\lead(J)$,
is the vector subspace of $S$ spanned by leading forms of all the elements of
$J$.
\end{definition}

For example, the inverse system $\langle x^3+y^2 \rangle = \Span\{x^3+y^2,x^2,x,y,1\}$ has
\begin{equation*}
  \lead(\langle x^3+y^2 \rangle) = \Span\{x^3,x^2,x,y,1\} = \langle x^3,y \rangle.
\end{equation*}
There is a tight connection between a system $J$ and $\lead(J)$.

\begin{proposition}\label{proposition: inverse system same hilbert function as leading forms}
The Hilbert functions of $J$ and $\lead(J)$ are equal.
\end{proposition}
\begin{proof}[sketch]
Let $f_1,f_2,\dotsc,f_s$ be a vector space basis for $\lead(J)$ consisting of homogeneous elements
and let $g_1,g_2,\dotsc,g_s \in J$ with $\lead(g_i) = f_i$.
One can show the $g_i$ are a basis for $J$.
Expressing the Hilbert functions of $J$ and $\lead(J)$ in terms of the $g_i$ and $f_i$ gives the result.
\qed
\end{proof}

The \emph{initial form} or \emph{lowest degree form} of a polynomial $g_i$ is its lowest degree homogeneous part.
The \emph{initial ideal} of an ideal $K$, denoted $\init(K)$,
is the ideal generated by the initial forms of all elements of $K$.

\begin{proposition}[{\cite[Proposition 3]{emsalem}}] \label{prop: apolar of lead = init of apolar}
Let $J$ be a finite-dimensional inverse system with ideal $J^\perp = I$.
We have $\lead(J)^\perp = \init(I)$.
In other words, $T/\lead(J)^{\perp}$ is the associated graded algebra of
$T/J^{\perp}$.
\end{proposition}

\begin{proof}
If $\Theta \in \init(I)$, then $\Theta = \init(\Psi)$, for some $\Psi \in I$.
To see that $\Theta \in \lead(J)^\perp$, let $F = \lead(G)$ for $G \in J$.
It follows that $\Theta \hook F$ is the highest degree part of $\Psi \hook G = 0$, so it
is zero. This shows that $\init(I) \subseteq \lead(J)^\perp$.
We have
\begin{multline*}
  \dim_{\kk} J = \dim_{\kk} \lead(J) = \dim_{\kk} T/\lead(J)^\perp \\
     \leq \dim_{\kk} T/\init(I) = \dim_{\kk} T/I = \dim_{\kk} J,
\end{multline*}
where the first equality is by Proposition \ref{proposition: inverse system same hilbert function as leading forms}
and the last is by Proposition \ref{proposition: inverse system same dimension as apolar algebra}.
This completes the proof.
\qed
\end{proof}

\begin{remark}
By Proposition~\ref{proposition: inverse system same hilbert function as leading forms}
and Proposition~\ref{prop: apolar of lead = init of apolar},
the Hilbert function of an inverse system $J$ is also the Hilbert function of a standard graded algebra,
namely the associated graded algebra of $T/J^\perp$. Hence, Macaulay's and
Gotzmann's theorems apply to these functions.
This enables us to prove that the only possible Hilbert functions $\vh$ of local ideals $I = J^\perp$ in $H_3^{11}$ with full embedding dimension $3$,
equivalently $\vh(1)=3$, are the ones listed
in Table~\ref{table: case proof sections}.
Since $\vh(2) \leq 6$,  we need to consider every possible value for $\vh(2)$, $1 \leq \vh(2) \leq 6$.
Also, $\sum \vh(i) = \dim_{\kk} T/I = 11$.
Finally, if $\vh(i) \leq 2$ for any $i \geq 2$, then $\vh$ is nonincreasing from the $i$th step onward,
by Corollary \ref{corollary: hilbert function nonincreasing}.
It is then easy to list the possible Hilbert functions and to check that all
of them are in Table~\ref{table: case proof sections}.
\end{remark}

\begin{proposition} [{\cite[\textsection C.2]{emsalem}}] \label{proposition: family of inverse systems}
Let $F(t) = \{f_1(t),f_2(t),\dotsc,f_s(t)\} \subset S[[t]]$
be a collection of polynomials in $S[[t]]$,
which we regard as polynomials in $S$ whose coefficients are continuous functions of a parameter $t$
in a neighborhood of $0$.
The family of apolar ideals $\{ F(t)^\perp \}$
satisfies $\lim_{t \to 0} F(t)^\perp \subseteq F(0)^\perp$.
If the inverse systems $\langle F(t) \rangle$ have the same Hilbert function for all $t$,
then we have $\lim_{t \to 0} F(t)^\perp = F(0)^\perp$ and $\{F(t)^\perp\}$ is a flat family.
\end{proposition}
\begin{proof}
If $\Theta \in \lim_{t \to 0} F(t)^\perp$,
write $\Theta = \Theta(0) = \lim_{t \to 0} \Theta(t)$ where $\Theta(t) \in F(t)^\perp$ for $t \neq 0$.
For each $t \neq 0$ we then have that $\Theta(t) \hook f_i(t) = 0$, for $i=1,\dotsc,s$.
By continuity, we also have that $\Theta(0) \hook f_i(0) = 0$.
This shows $\Theta \in F(0)^\perp$ and $\lim_{t \to 0} F(t)^\perp \subseteq F(0)^\perp$.
The equality of Hilbert functions implies equality of dimensions, so the ideals are equal.
\qed
\end{proof}

\begin{definition}\label{def: limit of inverse systems}
When $J_t = \langle f_1(t),f_2(t),\dotsc,f_s(t) \rangle$ is a parametrized family of inverse systems
generated by polynomials $f_i$ whose coefficients are continuous functions of $t$,
we will say $\lim_{t \to 0} J_t = J_0$ if and only if $\lim_{t \to 0} J_t^\perp = J_0^\perp$.
\end{definition}
\begin{example}\label{example: limit of inverse systems difference quotient}
Consider the families $W_1 = \{ \langle \ell^d,m^d \rangle \mid \ell,m \in S_1, \text{ independent}\}$
and $W_2 = \{ \langle \ell^d,\ell^{d-1}m \rangle \mid \ell,m \in S_1, \text{ independent}\}$.
Since the limit
\begin{equation*}
  \lim_{t \to 0} \frac{ (\ell + tm)^d - \ell^d }{ dt } = \ell^{d-1} m ,
\end{equation*}
we have, by Proposition \ref{proposition: family of inverse systems}, that
\begin{equation*}
  \lim_{t \to 0} \langle \ell^d, (\ell+tm)^d \rangle
    = \lim_{t \to 0} \left\langle \ell^d, \frac{(\ell+tm)^d-\ell^d}{dt} \right\rangle
    = \langle \ell^d, \ell^{d-1} m \rangle.
\end{equation*}
This is because every inverse system in each family has Hilbert function $(1,2,\dotsc,2)$.
This implies that $W_2$ is in the closure of $W_1$ in the Zariski topology.
\end{example}

\section{The Hilbert scheme of 11 points in 3-space}\label{section: Preparation}

In this section we,
use \sw{Macaulay2} to perform some computations that will be
needed later on and gather some general methods applicable to several
of the cases.

\subsection{\sw{Macaulay2} code examples}\label{section: Macaulay2 code examples}

To check if an ideal $I$ in $T=\kk[a,b,c]$
is smooth  we can run the following code.
This is one of the cases we check in the proof of Proposition \ref{case7}.
\begin{verbatim}
i1 : T = QQ[a,b,c]

i2 : I = ideal {b*c,a*b,a^2*c,a^3-c^2,b^5}

i3 : (dim I, degree I, degree Hom(I,T/I))

o3 = (0, 11, 33)

\end{verbatim}
These computations show that we have a zero-dimensional scheme of degree $11$ with tangent space dimension $33$.
If we now know that this is in the smoothable component, then it has to be a smooth point,
since we know that the smoothable component has dimension $3 \cdot 11 = 33$.
To check that this point is in the smoothable component, we construct a deformation.
We guess a candidate ideal $K$, then check that it satisfies the needed conditions.
\begin{verbatim}
i4 : R = T[t]

i5 : K = ideal {b*c,a*b,a^2*c,a^3-c^2,b^5+t*b^4}

i6 : assert (K:t == K)

i7 : minimalPrimes K

o7 = {ideal (c, a, t + b), ideal (c, b, a)}
\end{verbatim}
Here $K$ is an ideal in $\kk[a,b,c,t]$ whose special fiber (at $t=0$) is $I$.
To check that this is a flat family over $\kk[t]$, we appeal to \cite[Proposition III.9.7]{HarAG}
which implies that if the ideal $(K:t)$ equals $K$, then the family is flat in a neighbourhood of $0$.
The general fiber is supported at the two points $(0,-t,0),(0,0,0)$.
This shows the special fiber $I$ is cleavable, hence, by Lemma \ref{lemma: cleavable implies smoothable}
$I$ is also smoothable.

\subsection{Some general methods}

In this section we collect various results which we use in
Section~\ref{section: Cases}.

In our analysis of the irreducible components of some standard graded Hilbert scheme (and the fibers of $\pi_{\vh}$),
we will often consider the set of quadric generators $\{q_1,q_2,\dotsc, q_k\}$ of a homogeneous ideal $I\subset T$.
The following lemma describes the space of cubics $\langle q_1,q_2,\cdots, q_k\rangle\cdot T_1$ in the ideal generated
by these quadrics.

\begin{lemma}
\label{Quadrics via Gotzmann}
Let $T = \kk[\alpha_1,\alpha_2,\dotsc,\alpha_n]$ be the polynomial ring in $n$ variables.
Let $q_1,\dotsc ,q_k$ be linearly independent quadrics in $T$ where $2 \leq k \leq n$,
and let $I = (q_1,\dotsc,q_k)$.
Then $\dim I_3 \geq nk - \binom{k}{2}$,
with equality if and only if the $q_i$ share a common linear factor,
that is, $q_i = \ell \ell_i$ for some linear forms $\ell,\ell_1,\dotsc,\ell_k$.
\end{lemma}



\begin{proof}
Let $\vh$ be the Hilbert function of $T/I$.
The $2$-binomial expansion of $\vh(2)$ is given by
\(
\vh(2) = \binom{n+1}{2} - k = {\binom{n}{2} + \binom{n-k}{1}}.
\)
Thus,
\(
\vh(3) \leq \vh(2)^{\langle 2 \rangle} = {\binom{n+1}{3} + \binom{n-k+1}{2}},
\)
so
\begin{equation*}
  \dim I_3 = \dim T_3 - \vh(3) \geq \binom{n+2}{3} - \binom{n+1}{3} - \binom{n-k+1}{2} = nk - \binom{k}{2}.
\end{equation*}
Suppose that equality holds.
We will show that the $q_i$ share a linear factor.
By Gotzmann's Persistence Theorem, see Theorem \ref{Gotzmann's Persistence},
the equality
$\vh(3)=\vh(2)^{\langle  2 \rangle}$ implies that
$\vh(t+1) = \vh(t)^{\langle t \rangle}$ for all $t \geq 2$, which gives by induction
\begin{equation*}
  \vh(t) = \binom{n+t-2}{t} + \binom{n-k+t-2}{t-1} = \binom{n+t-2}{n-2} + \binom{n-k+t-2}{n-k-1}.
\end{equation*}
This shows that the projective scheme $V \subset \PP^{n-1}$ defined by $I$ has Hilbert polynomial of degree $n-2$
with leading coefficient $1/(n-2)!$.
By standard properties of Hilbert polynomials, see~\cite[Section I.7,
p.~52]{HarAG}, the scheme $V$ has codimension $1$ and degree $1$.
This means that $V$ consists of a reduced hyperplane $H$, possibly along with some lower-dimensional components.
Since each $q_i$ vanishes on $H$, they are all divisible by the equation $\ell$ of $H$.
\qed
\end{proof}

The following are generalizations of \cite[Proposition 4.3]{CEVV}.
\begin{lemma}\label{lemma: Hilbert scheme vector bundle over graded Hilbert scheme}
Fix $\vh = (1,\vh(1),\dotsc,\vh(t))$ with $h_i = \dim S_i$ for $i=1,\dotsc,t-2$.
The Hilbert scheme $H_n^{\vh}$ is a vector bundle of rank $\vh(t)(\dim S_{t-1} - \vh(t-1))$ over $\cH_n^{\vh}$.
In particular, the irreducible components of
$H_n^{\vh}$ are exactly the preimages of the irreducible components of $\cH_n^{\vh}$.
\end{lemma}
\begin{proof}
    A direct generalization of the proof for $t=3$ in \cite[Proposition 4.3]{CEVV}.
\qed
\end{proof}

\begin{lemma}\label{irreducibleFiber}
Fix $\vh = (1,\vh(1),\dotsc,\vh(t))$ with $\vh(i) = \dim T_i$ for $i=1,\dotsc,t-3$.
Every fiber of $\pi_{\vh}$ is isomorphic to an affine space; in particular, it is irreducible.
\end{lemma}
\begin{proof}
Let $I$ be a homogeneous ideal in $\cH^{\vh}_n$.
The fiber $\pi_{\vh}^{-1}(I)$ consists of ideals $I'$ with $\init(I') = I$
and with Hilbert function $\vh$.
Requiring that $\init(I')=I$ corresponds to adding higher degree terms to generators of $I$.
Requiring the Hilbert function of $T/I'$ to equal $\vh$ imposes conditions on the coefficients
of these higher degree terms.

Adding terms of degree greater than $t$ has no effect,
since these are already contained in $I$.
To any generator of degree $t-2$ or $t-1$, we can freely add terms of degree $t$
since they cannot change the Hilbert function.
To any degree $t-2$ generator $q_i$, we can add a term $a_i$ of degree $t-1$,
however, now there is something to check:
For any tuple of linear forms $\ell_1,\ell_2,\dotsc,\ell_r \in T_1$
such that $\ell_1 q_1 + \dotsb + \ell_r q_r = 0$,
we require that $\ell_1 a_1 + \dotsb + \ell_r a_r \in I'_t = I_t$.
These are all linear conditions on the coefficients of the $a_i$, hence, the solution space is an affine space.
Hence, the fiber at $I$ is isomorphic to $\A^k$ for some $k$.
\qed\end{proof}

\begin{remark}\label{rem:instead of Buchb. for 13421}
If there are only two generators $q_1$ and $q_2$ of degree $t-2$,
then there can be at most one (possibly trivial) linear condition on the forms $a_1$ and $a_2$ (as above).
Namely, if there are linear forms $\ell_1, \ell_2$ such that $\ell_1 q_1 + \ell_2 q_2 = 0$,
then these are uniquely determined up to a common scalar multiple,
and the condition $\ell_1 a_1 + \ell_2 a_2 \in I_t$ is sufficient for $\init(I') = I$.
\end{remark}

Going beyond the situation of Lemma~\ref{lemma: Hilbert scheme vector bundle over graded Hilbert scheme}, it is possible that the fibers of $\pi_{\vh}$ may be reducible.
To show that they are contained in the main component of the Hilbert scheme we would have to find
a smooth and smoothable point in each component of the fiber.
Unfortunately in general it is difficult to describe the fibers of $\pi_{\vh}$.
The following statement allows us in a handful of very special cases to avoid this difficulty.
\begin{lemma}\label{lemma: homogeneous lies in every component of fiber}
If $I \in \cH_n^{\vh}$, then $I$ lies in every irreducible component of $\pi_h^{-1}(I)$.
\end{lemma}
If the homogeneous ideal $I$ happens to be a smooth and smoothable point,
then the whole fiber is contained in the main component of the Hilbert scheme.
\begin{proof}
Let $I' \in \pi_{\vh}^{-1}(I)$, so that $I = \init(I')$.
The deformation of~\cite[Theorem~15.17]{eisenbud} gives a path in $\pi_{\vh}^{-1}(I)$ from $I'$ to $I$,
so $I$ lies in the irreducible component that contains $I'$.
\qed
\end{proof}

\subsection{Non-linear changes of coordinates}\label{sect: nonlinear}

We recall the technique of non-linear changes of coordinates
as in \cite{cjn13, EliasRossi_Analytic_isomorphisms} and \cite[Section 2.2]{JJlocArtGor}.
Assume we have a zero-dimensional quotient $A=T/I=\kk[\alpha,\beta,\gamma]/I$ supported at the origin.
The algebra $A$ can also be viewed as a quotient of the power series ring $R=\kk[[\alpha,\beta,\gamma]]$.
The power series ring has a much larger automorphism group than the polynomial ring.
Denote the maximal ideal of $R$ by $\mathfrak{m}$.
For any $\sigma_1,\sigma_2,\sigma_3 \in \mathfrak{m}$ whose images span $\mathfrak{m}/\mathfrak{m}^2$,
there is an automorphism $\phi$ of $R$ defined by $\phi(\alpha)=\sigma_1$, $\phi(\beta)=\sigma_2$, and $\phi(\gamma)=\sigma_3$.

Let $J=\langle f_1,f_2,\dotsc,f_r \rangle$ be the associated inverse system of $I$.
By \cite[Section 2.2]{JJlocArtGor}, the inverse system of $\phi^{-1}(I)$
is generated by $\phi^\vee(f_i)$ where $\phi^\vee$ is defined as follows.
Let $D_\alpha = \phi(\alpha) - \alpha$, $D_\beta = \phi(\beta) - \beta$, and $D_\gamma = \phi(\gamma) - \gamma$.
Then we have
\begin{equation*}
\phi^\vee(f) = \sum_{(k,m,n) \in \ZZ_{\geq 0}^3} \frac{x^ky^mz^n}{k!m!n!}
\cdot \left(D_\alpha^kD_\beta^mD_\gamma^n \hook f\right).
\end{equation*}

\begin{example}\label{example: nonlinear}
Let $J=\langle f \rangle$ for $f=x^4+y^4+g$ where $\deg g \leq 3$.
By subtracting multiples of $\alpha \hook f$ and $\alpha^2 \hook f$ from $f$,
we may assume the monomials $x^3$ and $x^2$ do not appear in $g$.
We will perform a non-linear change of coordinates so that there are
no monomials in $g$ divisible by $x^2$.
This will be needed in the proof of Lemma \ref{lem:case 13421,f->l^4+m^4}.

Let $B$ be the coefficient of $x^2y$ in $g$ and let $C$ be the coefficient of $x^2z$.
Let $\phi(\alpha)=\alpha$, $\phi(\beta) = \beta - \frac{B}{12}\alpha^2$, and $\phi(\gamma) = \gamma - \frac{C}{12} \alpha^2$.
We have $D_\alpha = 0$, $D_\beta = -\frac{B}{12}\alpha^2$, and $D_\gamma = -
\frac{C}{12} \alpha^2$, so
\begin{equation*}
\phi^\vee(x^4) = x^4 + y D_\beta(x^4) + zD_\gamma(x^4) + \dotsb =
x^4 -yBx^2-zCx^2 + \dotsb,
\end{equation*}
where we have omitted terms of degree less than $3$.
Similarly $\phi^\vee(y^4) = y^4$ and $\phi^\vee(g)$ is equal to $g$, modulo terms of degree less than $3$.
Also $\phi^\vee(f)$ will have no terms divisible by $x^2$.
\end{example}

\subsection{An explicit construction of flat families}

The section is adapted from \cite[Section 5]{cjn13}, where more general
results were proved for Gorenstein schemes.
Fix a zero-dimensional scheme $R$. In this section, under certain mild
assumptions on $R$, we construct a family with special fiber $R$ and general
fiber reducible, so that $R$ becomes cleavable.

    \begin{proposition}\label{ref:cleavinggeometric}
    Let $R \subset \A^n$ be a finite scheme supported at the
    origin. Let $C \subset \A^n$ be a smooth curve passing through
    the origin. Let $H = (x = 0)$ be a hyperplane intersecting $C$
    transversely.
    Let $r \geq 1$ be such that the ideal of intersection $R \cap C$ in $C$ is $(x^r)$
    and let $H^{r-1} = (x^{r-1} = 0)$ denote the thick hyperplane.
        If $R \subset C \cup H^{r-1}$ as schemes, then $R$ is cleavable.
    \end{proposition}
    \begin{proof}
        Since $R \cap C$ is cut out of $C$ by $x^r$, we can choose an $F \in I(R)$ whose image
        in $\kk[C]=\kk[A_n]/I(C)$ is $x^r$. Then we have $q := x^r-F \in I(C)$.  Now
        the image in $\kk[C]$ of any $i \in I(R)$  is $gx^r$, for some $g$. Write $i =
        g(x^r-q) + j$, for some $j$. We see that $j \in I(R) \cap I(C)$ which
        implies that $I(R) = (x^r-q) + I(R \cup C)$, hence, $R$ is cut out of $R \cup
        C$ by the equation $x^r-q$.
        There is a deformation of $R \subset R\cup C$ given by
        deforming this equation, namely
        \begin{equation}\label{eq:rayfamily}
            V(x^{r} - tx^{r-1} - q) \subset (R\cup C) \times A^1,
        \end{equation}
        with $t$ being the local parameter on $A^1$.

        To prove the flatness of the family \eqref{eq:rayfamily} it is enough to
        prove that every polynomial $f \in \kk[t]$ is not a zero-divisor in the coordinate ring of
        $V = V(x^r - tx^{r-1} - q)$.
        Suppose there is an $f\in \kk[t]$ and a function $g$ on $V$ such that $fg$ is zero.
        We will show that $g$ vanishes on $V\cap (C \times A^1)$ and on
        $V\cap (H^{r-1} \times A^1)$. Since $R \cup C
        \subset C\cup H^{r-1}$, this implies that $g$ vanishes on the whole of $V$,
        so that it is zero.

        First let us restrict to $C$, i.e. consider the family $V\cap (C \times
        A^1)$. It is given by the equation $x^r -tx^{r-1}$, thus, it is
        flat. Therefore, $f(t)$ is not a zero-divisor, hence, $g$ restricts to
        zero on $C \times A^1$.
        Next let us restrict to $H^{r-1}$, i.e. consider the family $V\cap
        (H^{r-1} \times A^1)$. It is given by the equation $x^r - q$,
        which does not involve $t$. Hence, this family is constant, thus,
        flat. Hence, $g$ restricts to zero on $H^{r-1} \times A^1$,
        which concludes the proof of flatness.
        The fiber of the family~\eqref{eq:rayfamily} over $t\neq 0$ is
        supported on at least two points: the origin and $(t, 0,  \ldots ,0)$,
        thus, reducible. Therefore, $R$ is cleavable.
    \qed
    \end{proof}

    \begin{corollary}\label{ref:algebraicversion}
        Let $R \subset A^n$ be a finite scheme supported at the
        origin. Let $I = I(R)$ be its ideal. Choose coordinates
        $\alpha_1,\alpha_2,
        \dotsc ,\alpha_n$ on $A^n$.
        Assume that $c$ is such that
        \(
            \alpha_1^c\cdot \alpha_j \in I(R)
        \)
        for all $j\neq 1$.
        Assume moreover that $\alpha_1^c\notin I + (\alpha_2,\alpha_3,\dotsc,\alpha_n)$.
        Then $R$ is cleavable.
    \end{corollary}
    \begin{proof}
        This follows from Proposition~\ref{ref:cleavinggeometric} above if we take $C =
        V(\alpha_2,\alpha_3,\dotsc,\alpha_n)$, $H = (\alpha_1)$. Then $r$ is defined by $R\cap C =
        (\alpha_1^r)$ and by assumption $r > c$, so that $R \subset C
        \cup H^{r-1}$.
    \qed
    \end{proof}

    \begin{corollary}\label{ref:unionoflineandpointissmoothable:cor}
        Suppose that $R \subset A^3$ is a scheme of length $11$. Let
        $I \subset \kk[\alpha, \beta, \gamma]$ be its ideal and suppose
        that $\alpha\beta, \alpha\gamma\in I$. Then the ideal $I$ is smoothable.
    \end{corollary}
    \begin{proof}
        If $R$ is reducible, it is smoothable because all its
        components are. Suppose $R$ is irreducible supported at the origin.
        If any order one element lies in $I$, then after a non-linear
        coordinate change $R$
        is contained in an $A^2$ and so is smoothable. If no order
        one element lies in $I$, then Corollary~\ref{ref:algebraicversion}
        applied to $c = 1$ implies that $R$ is cleavable. Therefore, it is
        smoothable by Lemma~\ref{lemma: cleavable implies smoothable}.
        \qed
    \end{proof}

\section{Proof of main theorem}\label{section: Cases}

In this section, we prove Theorem~\ref{ref:main:theorem} by proving for each possible Hilbert function
that algebras with that Hilbert function are smoothable.
For reference, Table~\ref{table: case proof sections} shows in which section each Hilbert function is treated.

In this section we fix $n=3$, $S = \kk[x,y,z]$, and $T = \kk[\alpha,\beta,\gamma]$.

\subsection{Cases with long tails of ones}\label{sec:tailsofones}

\begin{proposition}
Let $\vh$ be one of these Hilbert functions:
$(1,3,1,1,1,1,1,1,1)$,
$(1,3,2,1,1,1,1,1)$,
$(1,3,2,2,1,1,1)$,
$(1,3,3,1,1,1,1)$,
$(1,3,4,1,1,1)$,
$(1,3,5,1,1)$.
Then we have $H_3^{\vh} \subset R_3^{11}$.
\end{proposition}
\begin{proof}
Let $I \in H_3^{\vh}$ be an ideal with Hilbert function $\vh$.
By Proposition~\ref{ref:tailsofones:prop}, the ideal $I$ is cleavable.
By Lemma~\ref{lemma: cleavable implies smoothable}, the ideal $I$ is smoothable.
So $I \in R_3^{11}$.
\qed
\end{proof}

\begin{proposition}\label{ref:tailsofones:prop}
    Let $R = \Spec A \subset A^n$ be an irreducible subscheme and
    $\vh$ be the Hilbert function of the local algebra $A$.
    Suppose $\vh = (1, \vh(1),  \ldots , \vh(c), 1,  \ldots , 1)$ with at least
    $c$ trailing ones,
    that is, letting $s$ be the greatest value such that $\vh(s) \neq 0$,
    we assume that $\vh(k) = 1$ for $c+1 \leq k \leq s$, and $s \geq 2c$.
    Then the scheme $R$ is cleavable.
\end{proposition}
The proof follows the Gorenstein case of \cite[Example 5.15]{cjn13}.
\begin{proof}
    Let $I$ be the ideal of $R$ and let $J$ be the inverse system of $I$.
    Consider a minimal generating set of $J$.
    It has a unique generator $f$ of degree $s$.
    As explained in Section~\ref{sect: nonlinear}, we can perform a non-linear coordinate change to assume that
    \(
      f = x_1^s + g,
    \)
    for some
    $g$ such that $\alpha_1^c \hook g = 0$. All other  generators of $J$ are of degree
    at most $c$. By subtracting some partials of $f$, we may assume that they
    are also annihilated by $\alpha_1^c$. Thus, $\alpha_1^c \alpha_j$ lies in $I$ for all $j \neq 1$.

    It remains to check that $\alpha_1^c\notin I + (\alpha_2, \alpha_3, \ldots, \alpha_n)$.
    Take any $q\in (\alpha_2, \alpha_3,\dotsc ,\alpha_n)$.
    Then $(\alpha_1^c - q) \hook f = \frac{s!}{(s-c)!} x_1^{s-c} - q \hook g$.
    We claim this is nonzero.
    Note that $s-c \geq c$ by assumption on the number of trailing ones.
    Therefore, $\alpha_1^{s-c}$ annihilates $g$.
    So
    \begin{equation*}
      \alpha_1^{s-c} \hook \left( \frac{s!}{(s-c)!} x_1^{s-c} - q \hook g \right) = s! x_1^0 - q \hook (\alpha_1^{s-c} \hook g) = s! \neq 0 .
    \end{equation*}
    This shows $(\alpha_1^c - q) \hook f \neq 0$, as claimed,
    so $\alpha_1^c - q \notin I$.
    Therefore, $\alpha_1^c \notin I + (\alpha_2,\alpha_3,\dotsc,\alpha_n)$.
    Thus, by Corollary~\ref{ref:algebraicversion} the subscheme $R$ is
    cleavable.
\qed
\end{proof}

\subsection{Cases with short Hilbert functions}\label{sec:1def}

For the three cases $\vh=(1,3,3,4)$, $\vh=(1,3,4,3)$, and $\vh=(1,3,5,2)$,
the analysis of the irreducible components of their standard graded Hilbert schemes completely determines the corresponding strata
in the (not graded) Hilbert scheme $H_3^{\vh}$.
Explicitly, in each of these cases $H_3^{\vh}$ is a vector bundle over $\cH_3^{\vh}$
by Lemma \ref{lemma: Hilbert scheme vector bundle over graded Hilbert scheme},
so the irreducible components of $H_3^{\vh}$ are exactly the preimages of the irreducible components of $\cH_3^{\vh}$.

In each of the three cases, we will first cover $\cH_3^{\vh}$
by a collection of irreducible sets (which are not necessarily components)
and produce a smooth and smoothable ideal for each set.
By Lemma \ref{lemma: Hilbert scheme vector bundle over graded Hilbert scheme}
and Lemma~\ref{lemma: smooth and smoothable is enough},
this is enough to guarantee that all algebras in $H_3^{\vh}$ are smoothable.

\begin{proposition}\label{proposition: 1334}
Let $\vh = (1,3,3,4)$.
Then $H_3^{\vh} \subset R_3^{11}$.
\end{proposition}
\begin{proof}
Let $I \subset T$, $I \in \cH_3^{\vh}$ be a homogeneous ideal
such that $A = T/I$ has Hilbert function $\vh$.
Then $\dim I_2 = \dim T_2 - \vh(2) = 3$.
Let $I' = (I_2)$ be the ideal generated by the quadrics in $I$.
By Lemma~\ref{Quadrics via Gotzmann}, $\dim I'_3 \geq 3 \cdot 3 - \binom{3}{2} = 6$,
but $\dim I'_3 \leq \dim I_3 = \dim T_3 - \vh(3) = 6$.
So $I_3 = I'_3$, equality holds in the dimension bound,
and by Lemma~\ref{Quadrics via Gotzmann}, the quadrics in $I_2$ must share a common linear factor $\ell$.

Then $I_2$ is spanned by $\ell \alpha$, $\ell \beta$, $\ell \gamma$.
That is, the standard graded Hilbert scheme $\cH_3^{\vh}$ is parametrized by the line $\ell$.
It is, therefore, isomorphic to the Grassmannian $\Gr(1,3)\cong \PP^2$ and, hence, irreducible.
By Lemma \ref{lemma: Hilbert scheme vector bundle over graded Hilbert scheme},
$H_3^{\vh}$ is also irreducible.

It is sufficient to find one smooth and smoothable point in $H_3^{\vh}$.
Consider the ideal
$L=(\alpha\beta,\alpha\gamma,\alpha^2+\beta^3,\beta^2\gamma^2,\beta\gamma^3,\gamma^4)$.
It is smoothable by Corollary~\ref{ref:unionoflineandpointissmoothable:cor}
and we check computationally that $L$ is smooth.
\qed
\end{proof}

\begin{proposition}
Let $\vh = (1,3,4,3)$.
Then $H_3^{\vh} \subset R_3^{11}$.
\end{proposition}
\begin{proof}
The standard graded Hilbert scheme $\cH_3^{\vh}$ is a union of two irreducible sets.
We will provide a smooth and smoothable point in each of them.

Let $I \subset T$, $I \in \cH_3^{\vh}$ be a homogeneous ideal
such that $A = T/I$ has Hilbert function $\vh$.
Then $\dim I_2 = 2$.
By Lemma~\ref{Quadrics via Gotzmann},
the space of cubics generated by the quadrics in $I_2$ can have dimension either $6$ or $5$,
and the latter occurs exactly when the quadrics share a linear factor.
Let $P \subset \cH_3^{\vh}$ be the set of ideals $I$ whose quadrics generate a $6$-dimensional space of cubics
and let $Q \subset \cH_3^{\vh}$ be the set of ideals $I$ whose quadrics generate a $5$-dimensional space of cubics.
Then $\cH_3^{\vh} = P \cup Q$.
We claim that each of $P$ and $Q$ is irreducible.

The subset $P$ is parametrized by pairs of spaces $(K,M)$,
where $K$ is a $2$-dimensional subspace of $T_2$, 
not of the form $\Span\{\ell \cdot \ell_1, \ell \cdot \ell_2\}$,
and $M$ is a $7$-dimensional subspace of $T_3$ that contains $K\cdot T_1$,
equivalently a line in $T_3/K \cdot T_1$.
Thus, $P$ is realized as a projective bundle with fiber $\PP(T_3/K \cdot T_1)$
over an open subset of $\Gr(2,T_2)$.
In particular, $P$ is irreducible.

In the subset $Q$, the quadrics $q_1,q_2$ that span $I_2$
have the form $q_1 = \ell \cdot \ell_1$ and $q_2 = \ell \cdot \ell_2$ for some lines $\ell, \ell_1, \ell_2$.
This component is parametrized by a triple $(\ell,L,N)$, where $\ell \in T_1$ is the common line,
$L=(\ell_1,\ell_2)\subset T_1$ is the space spanned by the other two lines,
and $N$ is a $7$-dimensional space of $T_3$ that contains the $5$-dimensional space $\ell \cdot L\cdot T_1$.
So $Q$ is isomorphic to a Grassmannian bundle with fiber
$\Gr(7-5,T_3/\ell \cdot L \cdot T_1)$,
over a base $\Gr(1,T_1)\times\Gr(2,T_1)$; it is, therefore, irreducible.

Now $H_3^{\vh} = \pi_{\vh}^{-1}(P) \cup \pi_{\vh}^{-1}(Q)$,
and by Lemma \ref{lemma: Hilbert scheme vector bundle over graded Hilbert scheme} these are irreducible sets as well.
To complete this case, we provide a smooth and smoothable ideal for each set.
The ideal $I=(\alpha^2,\beta^2,\gamma^3,\alpha\beta\gamma^2)$ lies in $P$ and,
hence, also in $\pi_{\vh}^{-1}(P)$.
It is monomial, hence, smoothable by \cite[Proposition~4.15]{CEVV}
and it is easy to check computationally that it is a smooth point.
For $\pi_{\vh}^{-1}(Q)$ let 
\(
  I=(\alpha\beta,\alpha\gamma,\alpha^3+\gamma^3,\beta\gamma^2,\beta^3\gamma,\beta^4).
\)
Then $I$ is smoothable by
Corollary~\ref{ref:unionoflineandpointissmoothable:cor}
and once again a smooth point.
\qed
\end{proof}

\begin{proposition}\label{proposition: 1352}
Let $\vh=(1,3,5,2)$.
Then $H_3^{\vh} \subset R_3^{11}$.
\end{proposition}
\begin{proof}
Let $I \subset T$, $I \in \cH_3^{\vh}$ be a homogeneous ideal
such that $A = T/I$ has Hilbert function $\vh$.
Then $\dim I_2 = 1$ and $\dim I_3 = \dim T_3 - \vh(3) = 8$.
The standard graded Hilbert scheme $\cH^3_{\vh}$ is parametrized by pairs $(L,M)$,
where $L$ is some $1$-dimensional subspace of $T_2$
and $M$ is an $8$-dimensional subspace of $T_3$ that contains the $3$-dimensional subspace $L\cdot T_1$. 
This parametrization realizes an isomorphism of
$\cH_3^{\vh}$ to a Grassmannian bundle with base $\PP T_2$ and fiber $\Gr(8-3,T_3/L \cdot T_1)$,
proving that $\cH_3^{\vh}$ is irreducible.
By Lemma \ref{lemma: Hilbert scheme vector bundle over graded Hilbert scheme},
$H_3^{\vh}$ is irreducible as well.

Now let
\(
  I = (\alpha\beta,\alpha^3,\beta^3,\gamma^3,\alpha\gamma^2,\alpha^2\gamma+\beta\gamma^2).
\)
One can check that $I \in H_3^{\vh}$.
Since $\alpha\beta^2, \alpha\gamma^2 \in I$ and $\alpha^2 \notin I+(\beta,\gamma)$,
Corollary~\ref{ref:algebraicversion} with $c=2$ implies $I$ is smoothable.
Finally one can check computationally that $I$ is a smooth point.
\qed
\end{proof}

\subsection[Case (1,3,4,2,1)]{Case $\vh=(1,3,4,2,1)$}\label{sec:13421}

\begin{proposition}
Let $\vh=(1,3,4,2,1)$.
Then $H_3^{\vh} \subset R_3^{11}$.
\end{proposition}

\begin{proof}
Let $I \in \cH_3^{\vh}$ be a homogeneous ideal with inverse system $J$.
Let $f \in J_4$ and let $\vh_f$ be the Hilbert function of $\langle f \rangle$.
Since $\langle f \rangle \subset J$ we have $\vh_f \leq \vh$.
By Proposition~\ref{ref:symmetry} and Macaulay's bound (Theorem~\ref{Macaulay's bound}),
$\vh_f$ must be $(1,2,3,2,1)$, $(1,2,2,2,1)$, or $(1,1,1,1,1)$.

If $\vh_f = (1,2,3,2,1)$
see Lemma \ref{lem:case 13421,f->12321}.
If $\vh_f = (1,2,2,2,1)$ then by Remark \ref{remark: *11 *221} we can choose coordinates so that
$f = x^4+y^4$ or $f = x^3 y$.
For $f=x^4+y^4$ see Lemma \ref{lem:case 13421,f->l^4+m^4} and for $f=x^3y$ see Lemma \ref{lem:case 13421,f->l^3m}.
If $\vh_f = (1,1,1,1,1)$ see Lemma \ref{lem:case 13421,f->11111}.
\qed
\end{proof}

\begin{lemma}
\label{lem:case 13421,f->12321}
Let $\vh = (1,3,4,2,1)$ and let $I \in \cH^{\vh}_3$ be a homogeneous ideal
with inverse system $J$.
Suppose that the degree $4$ generator $f$ of $J$ is such that the Hilbert function
of $\langle f \rangle$ is $(1,2,3,2,1)$.
Then $\pi_{\vh}^{-1}(I) \subset R^{11}_3$.
\end{lemma}

\begin{proof}
Let $I' \in \pi_{\vh}^{-1}(I)$ with inverse system $J'$.
Let $F$ be the degree $4$ generator of $J'$, so that $f$ is the leading form of $F$.
We will construct a family $J'_t$ so that $J'_1 = J'$ and $J'_0$ is $\langle x^2y^2,z^2 \rangle$,
$\langle x^2y^2,zx \rangle$, or $\langle x^2y^2,z(x+y) \rangle$.
First change coordinates so that $f \in \kk[x,y]$.
Then $x^2,xy,y^2 \in \langle f \rangle_2 \subset J'$,
so $J'_{\leq 2}$ is spanned by $\{x^2,xy,y^2,Q,S_{\leq 1}\}$ for a quadratic form $Q \in \kk[x,y,z]$.
Write $Q = cxz + dyz + ez^2$.
If $e \neq 0$ then changing coordinates by replacing $z$ with a suitable linear combination of $x,y,z$ to complete the square
eliminates the $xz$ and $yz$ terms and takes $Q$ to $z^2$ modulo $x^2,xy,y^2$.
So either $J' = \langle F,z^2 \rangle$ or $J' = \langle F,z(cx+dy) \rangle$.

Write $F = f + g$, $\deg g \leq 3$.
By well-known facts about binary forms (see for example, \cite[Theorem 1.43]{iakanev}),
we have $f = \ell_1^4 + \ell_2^4 + \ell_3^4$ for some nonproportional linear forms $\ell_i \in \kk[x,y]$.
Observe that $18 x^2 y^2 = (x+y)^4 + \omega (x+\omega y)^4 + \omega^2(x+\omega^2 y)^4$ where $\omega$ is a cube root of unity.
We change coordinates in $\kk[x,y]$ so that $\ell_1 = x+y$ and $\ell_2 = \omega^{1/4}(x+\omega y)$.
Let $f_t = \ell_1^4 + \ell_2^4 + (t \ell_3 + (1-t)\omega^{1/2}(x+\omega^2 y))^4$, $F_t = f_t + tg$,
and $J'_t = \langle F_t,Q \rangle$.

It is easy to check that $F_1 = F$, $F_0 = 18 x^2 y^2$, and for all but finitely many $t$,
$\langle f_t \rangle$ has Hilbert function $(1,2,3,2,1)$ and $J'_t$ has Hilbert function $(1,3,4,2,1)$.
Then $\lim J'_t = J'_0 = \langle 18 x^2 y^2, Q \rangle = \langle x^2 y^2, Q \rangle$, as in Definition \ref{def: limit of inverse systems}.
Rescaling $x$ and $y$ and interchanging if necessary, $Q$ is one of $z^2$, $zx$, or $z(x+y)$.
Now $\langle x^2y^2,z^2 \rangle^\perp$ and $\langle x^2y^2,zx \rangle^\perp$ are monomial ideals, hence, smoothable.
The family $(\gamma^2,\alpha\gamma-\beta\gamma,\beta^2\gamma,\beta^3,\alpha^3 + t \alpha^2)$
shows that 
$\langle x^2y^2,z(x+y) \rangle^\perp = (\gamma^2,\alpha\gamma-\beta\gamma,\beta^2\gamma,\beta^3,\alpha^3)$ is smoothable.
So all three points are smoothable and it is easy to check that each one is a smooth point.
Hence, the irreducible (one-dimensional) family $\{(J'_t)^\perp\} \subset R^{11}_3$, in particular $I' = (J'_1)^\perp \in R^{11}_3$.
\qed
\end{proof}

\begin{lemma}
\label{lem:case 13421,f->l^4+m^4}
Let $\vh=(1,3,4,2,1)$ and let $I \in \cH^{\vh}_3$ be a homogeneous ideal with inverse system $J$.
Suppose that the degree $4$ generator of $J$ is of the form $\ell^4+m^4$ for some independent linear forms
$\ell,m \in S_1$.
Then $\pi_{\vh}^{-1}(I) \subset R^{11}_3$.
\end{lemma}

\begin{proof}
    Assume $\ell=x,m=y$.
    Let $I' \in \pi_{\vh}^{-1}(I)$ with inverse system $J'$.
    We will apply Corollary~\ref{ref:algebraicversion}.
    Consider the degree four generator $F = x^4 + y^4 + g\in J'$, where $\deg g\leq 3$.
    Since $x^2 \in J$ we can subtract the $x^2$ term out of $g$. Then the
    only terms of $g$ divisible by $x^2$ are possibly $x^3$, $x^2y$, $x^2z$.
    After a non-linear coordinate change as in Example~\ref{example: nonlinear} we
    may assume that there are no such terms. Then $\alpha^2\hook F = 12x^2$, so $\alpha^2\not\in
    F^{\perp} + (\beta, \gamma)$. Moreover  $\alpha^2\beta$ and
    $\alpha^2\gamma$ annihilate $F$ and so its partials, hence, lie in $I'$.
    Therefore, the assumptions of Corollary~\ref{ref:algebraicversion} for
    $c=2$ are satisfied and $I'$ is cleavable. By Lemma~\ref{lemma:
    cleavable implies smoothable}, it is smoothable.
    \qed
\end{proof}

\begin{lemma}
\label{lem:case 13421,f->l^3m}
Let $\vh=(1,3,4,2,1)$ and let $I \in \cH^{\vh}_3$ be a homogeneous ideal with inverse system $J$.
Suppose that the degree $4$ generator of $J$ is of the form $\ell^3 m$ for some independent linear forms
$\ell,m \in S_1$.
Then $\pi_{\vh}^{-1}(I) \subset R^{11}_3$.
\end{lemma}

\begin{proof}
Assume $\ell=x$, $m=y$, so that $J = \langle x^3 y, Q_1, Q_2 \rangle$ for some quadratic forms $Q_1, Q_2$.
Let $I' \in \pi_{\vh}^{-1}(I)$ with inverse system $J'$.
We will show $I'$ is smoothable by writing it as a limit of smoothable points.
Note, $J' = \langle x^3 y + g_3 + g_2, Q_1, Q_2 \rangle$ where $g_i$ is a form of degree $i$ for $i=2,3$.
We introduce a parameter $t$ and let $y_t = x+ty$.
Observe that $\lim_{t\to 0} (y_t^4-x^4)/4t = x^3 y$.
For general $t$ we will define a form $g_3(t)$ so that $J'_t = \langle (y_t^4-x^4)/4t + g_3(t) + g_2, Q_1, Q_2 \rangle \to J'$
in the sense of Definition \ref{def: limit of inverse systems}.

To define $g_3(t)$, first note that $\gamma \hook g_3 \in J_2 = \Span\{x^2,xy,Q_1,Q_2\}$.
For $i=1,2$ let $Q^{\sharp}_i = \int Q_i \, dz$ be a homogeneous form of degree $3$ so that $\gamma \hook Q^{\sharp}_i = Q_i$.
Write $g_3 = a x^2z + b xyz + c Q^{\sharp}_1 + d Q^{\sharp}_2 + e(x,y)$
for some scalars $a,b,c,d$ and a $3$-form $e$.
Now we define
$g_3(t) = a x^2z + (b/2t)(y_t^2 - x^2)z + c Q^{\sharp}_1 + d Q^{\sharp}_2 + e(x,y)$.

Now $\gamma \hook g_3(t) \in \Span\{x^2,y_t^2,Q_1,Q_2\}$,
hence, $\lead(J'_t) = \langle y_t^4-x^4, \allowbreak Q_1, \allowbreak Q_2 \rangle$.
Since $\dim J_2 = 4$ we have $xy\not\in \Span\{x^2,Q_1,Q_2\}$.
Since $xy = \lim_{t\to 0} (y_t^2 - x^2)/(2t)$ we also have $y_t^2\not\in
\Span\{x^2,Q_1,Q_2\}$ for general $t$.
For such $t$ the space $(\lead(J'_t))_2$ has dimension $4$,
which means that $\lead(J'_t)$ and $J'_t$ have Hilbert function $\vh$.
Also $\lim_{t\to 0} g_3(t) = g_3$.
Therefore, $\lim_{t \to 0} J'_t = J'$, as desired.
By Lemma \ref{lem:case 13421,f->l^4+m^4}, each $(J'_t)^\perp$ with $t\neq 0$
 is smoothable, which implies $I' = \lim (J'_t)^\perp$ is smoothable as well.
\qed
\end{proof}

\begin{lemma}
\label{lem:case 13421,f->11111}
Let $\vh=(1,3,4,2,1)$ and let $I \in \cH^{\vh}_3$ be a homogeneous ideal
with inverse system $J$.
Suppose that the degree $4$ generator $f$ of $J$ is such that the Hilbert function
of $\langle f \rangle$ is $(1,1,1,1,1)$.
Then $\pi_{\vh}^{-1}(I) \subset R^{11}_3$.
\end{lemma}

\begin{proof}
By Remark \ref{remark: *11 *221}, we can choose coordinates so that $f = z^4$.
Let $V \subset \cH^{\vh}_3$ be the set of ideals $I$ satisfying the hypothesis, that is,
\(
  V = \{ I \in \cH^{\vh}_3 \mid I \subset (z^4)^\perp \}.
\)
For $I \in V$, $\dim I_2 = 2$ and $\dim I_3 = 8$.
By Lemma~\ref{Quadrics via Gotzmann}, $\dim T_1 \cdot I_2$ is either $5$ or $6$.
Let $V_1 \subset V$ be the set of $I$ such that $\dim T_1 \cdot I_2 = 6$,
equivalently the quadrics in $I_2$ have no common factor.
Let $V_2 \subset V$ be the set of $I$ such that $\dim T_1 \cdot I_2 = 5$
and $I_2 = \Span\{\ell \ell_1, \ell \ell_2\}$ for some linear forms $\ell, \ell_1, \ell_2$
such that $\Span\{\ell_1,\ell_2\} \subseteq z^\perp = \Span\{\alpha,\beta\}$ (necessarily equality must hold).
And let $V_3 \subset V$ be the remainder, the set of $I$ such that $\dim T_1 \cdot I_2 = 5$
and $I_2 = \Span\{\ell \ell_1, \ell \ell_2\}$ for some linear forms $\ell, \ell_1, \ell_2$
such that $\Span\{\ell_1,\ell_2\} \not\subset z^\perp$.
We will show that each $V_i$ and each $\pi_{\vh}^{-1}(V_i)$ is irreducible,
and give a smooth and smoothable point on each $\pi_{\vh}^{-1}(V_i)$.

First, every ideal $I \in V$ is determined by $(I_2,I_3)$.
Suppose $I \in V_1$.
The subspace $I_2 \subset (z^2)^\perp$ is parametrized by an open subset of $\Gr(2,(z^2)^\perp) = \Gr(2,5)$.
And then $I_3 \subset (z^3)^\perp$ is such that $T_1 \cdot I_2 \subset I_3$.
The quotient $I_3/T_1 \cdot I_2$ is a $2$-dimensional subspace of $(z^3)^\perp/T_1 \cdot I_2$.
So for each choice of $I_2$, $I_3$ may be chosen from $\Gr(8-6,(z^3)^\perp/T_1 \cdot I_2) = \Gr(2,3)$.
This shows $V_1$ is a Grassmannian bundle over an open subset of a Grassmannian,
in particular irreducible.
Let $I_2 = \Span\{q_1,q_2\}$.
Since $I\in V_1$, there are no lines $\ell_1,\ell_2$ such that $\ell_1 q_1 + \ell_2 q_2 = 0$.
By Lemma \ref{irreducibleFiber} and Remark \ref{rem:instead of Buchb. for 13421},
the fiber $\pi_{\vh}^{-1}(I)$ is a certain product of affine spaces.
Explicitly it is $T_3^2 \times T_4^4$, corresponding to cubic terms that may be added to the quadric generators of $I$
and quartic terms that may be added to the quadric and cubic generators of $I$.
This makes $\pi_{\vh}^{-1}(V_1)$ a (trivial!) vector bundle over $V_1$, hence, irreducible.
A smooth and smoothable point in $\pi_{\vh}^{-1}(V_1)$ is given by
$\langle yz, x^2 y, z^4 \rangle^\perp = (\alpha\gamma,\beta^2,\beta\gamma^2,\alpha^3,\gamma^5)$.
It is smoothable because it is a monomial ideal and we check computationally that it is a smooth point.
This shows that $\pi_{\vh}^{-1}(V_1) \subset R^{11}_3$.

If $I \in V_2$ then $I_2 = \ell \cdot \Span\{\alpha,\beta\}$ for some linear form $\ell$,
so $I_2$ is determined by the choice of $[\ell] \in \PP T_1$.
As before, for each choice of $\ell$, $I_3$ may be chosen from $\Gr(8-5,(z^3)^\perp/T_1 \cdot I_2) = \Gr(3,4)$.
Again this makes $V_2$ a Grassmannian bundle over an irreducible base, so $V_2$ is irreducible.
By Remark \ref{rem:instead of Buchb. for 13421},
$\pi_{\vh}^{-1}(V_2)$ is a trivial subbundle of a trivial vector bundle over $V_2$,
namely $\pi_{\vh}^{-1}(V_2) \subset V_2 \times (T_3^2 \times T_4^5)$
is defined by $\beta a_1 - \alpha a_2 \in I_4 = (z^4)^\perp$, where $a_1,a_2$ are the cubic terms added to the quadric generators
$\ell \alpha, \ell \beta$.
Hence, $\pi_{\vh}^{-1}(V_2)$ is irreducible.
A smooth and smoothable point in this set is given by the limit of the flat family
$(\alpha\gamma,\beta\gamma,\beta^3+\gamma^4,\alpha^3-t\cdot\alpha^2,\alpha^2\beta)$.

If $I \in V_3$ then, writing $I_2 = \Span\{\ell \ell_1, \ell \ell_2\}$, we must have $\ell \hook z = 0$,
since for at least one of $i=1,2$ we have $\ell_i \hook z \neq 0$, but $\ell \ell_i \hook z^2 = 0$.
Now $\ell$ may be chosen from $z^\perp$ and $\Span\{\ell_1, \ell_2\}$ may be chosen to be any $2$-dimensional subspace
of $T_1$ other than $z^\perp$.
So the choice of $I_2$ is parametrized by an open subset of $\PP(z^\perp) \times \Gr(2,T_1)$.
Once again, for each choice of $I_2$, $I_3$ may be chosen from the Grassmannian $\Gr(8-5,(z^3)^\perp/T_1 \cdot I_2) = \Gr(3,4)$.
Hence, $V_3$ is a Grassmannian bundle over an irreducible base, in particular irreducible.
By Remark \ref{rem:instead of Buchb. for 13421},
$\pi_{\vh}^{-1}(V_3)$ is a (nontrivial) subbundle of a trivial vector bundle over $V_2$,
namely $\pi_{\vh}^{-1}(V_3) \subset V_3 \times (T_3^2 \times T_4^5)$
is defined by $\ell_2 a_1 - \ell_1 a_2 \in I_4 = (z^4)^\perp$ where, as before, $a_1,a_2$ are the cubic terms
added to the quadric generators $\ell \ell_1, \ell \ell_2$.
Hence, $\pi_{\vh}^{-1}(V_3)$ is irreducible.
The ideal
$(\beta^2,\beta\gamma,\alpha^3,\alpha^2\gamma,\alpha\gamma^2,\gamma^5) \in V_3$
is smoothable because it is monomial and we check computationally that it is smooth.
\qed
\end{proof}

\subsection[Case (1,3,2,2,2,1)]{Case $\vh=(1,3,2,2,2,1)$}\label{sec:132221}

\begin{lemma}\label{lemma: 221 irreducible}
Let $\vh = (1,\vh(1),\dotsc,\vh(k),2,\dotsc,2,1)$
such that $\vh(i) = \dim S_i$ for all $i \leq k$, then has at least two $2$s and  a $1$  in the last position.
Then the standard graded Hilbert scheme $\cH_n^{\vh}$ is irreducible.
Each ideal $I \in \cH_n^{\vh}$ is the apolar ideal $J^\perp$ of an inverse system $J$
of one of the following forms: $\langle \ell^d + m^d, S_k \rangle$, $\langle \ell^{d-1}m, S_k \rangle$,
$\langle \ell^d, m^{d-1}, S_k \rangle$, $\langle \ell^d, \ell^{d-2}m, S_k \rangle$
for some linear forms $\ell, m$.
\end{lemma}
\begin{proof}
Say the last $1$ is in degree $d$,
let $J$ be a homogeneous inverse system with Hilbert function $\vh$, and let $f \in J$ be the $d$-form that appears.
Either $\langle f \rangle$ has Hilbert function $(\dotsc,2,2,1)$ or $(\dotsc,1,1,1)$.
In the first case $f = \ell^d + m^d$ or $f = \ell^{d-1} m$,
and $J$ is generated by $f$ together with $S_k$.
The second type is a limit of the first type,
similarly to Example \ref{example: limit of inverse systems difference quotient}.

In the second case $f = \ell^d$ and there is a generator $g$ of degree $d-1$.
Note $g$ has at most $2$ first derivatives since $\langle g \rangle_{d-2} \subseteq J_{d-2}$.
So $\langle g \rangle$ has Hilbert function $(\dotsc,2,1,0)$ or $(\dotsc,1,1,0)$.
If it is $(\dotsc,1,1,0)$ then $g = m^{d-1}$ for a linear form $m$ independent from $\ell$.
If the Hilbert function of $g$ is $(\dotsc,2,1,0)$ then $\ell^{d-2} \in \langle g \rangle_{d-2}$,
so $g = \ell^{d-2} m$ for a linear form $m$ independent from $\ell$.

So either $g = \ell^{d-2} m$ or $g = m^{d-1}$.
Correspondingly, either $J = \langle \ell^d, \ell^{d-2} m, S_k \rangle$ or $J = \langle \ell^d, m^{d-1}, S_k \rangle$.
Both of these can be obtained as limits of inverse systems of the
first two forms in appropriate ways,
using Proposition \ref{proposition: family of inverse systems}.
Explicitly, $\langle \ell^d, \ell^{d-2}m, S_k \rangle = \lim_{t \to 0} \langle \ell^{d-1}(\ell + t m), S_k \rangle$
and $\langle \ell^d, m^{d-1}, S_k \rangle = \lim_{t \to 0} \langle \ell^d + t m^d, S_k \rangle$.
\qed
\end{proof}

\begin{proposition}\label{proposition: 132221}
Let $\vh = (1,3,2,2,2,1)$. Then $H_3^{\vh} \subset R_3^{11}$.
\end{proposition}
\begin{proof}
By Lemma \ref{lemma: 221 irreducible},
every homogeneous ideal in $\cH_3^{\vh}$ is the apolar ideal
of an inverse system which is isomorphic to one of the following:
$J_1 = \langle x^5+y^5,z \rangle$, $J_2 = \langle x^4 y, z \rangle$,
$J_3 = \langle x^5, y^4, z \rangle$, or $J_4 = \langle x^5, x^3 y, z \rangle$.
We may dispose of the first two cases easily.
We compute $I_2 = J_2^\perp = (\alpha^5,\beta^2,\alpha\gamma,\beta\gamma,\gamma^2)$.
Then $I_2$ is smoothable because it is a monomial ideal and one can easily check computationally
that it is a smooth point.
By Lemma \ref{lemma: homogeneous lies in every component of fiber},
the smooth and smoothable point $I_2$ lies in every component of the fiber $\pi_{\vh}^{-1}(I_2)$,
which shows that each irreducible component of the fiber is contained in $R_3^{11}$.

Similarly, $I_1 = J_1^\perp = 
(\alpha^5-\beta^5,\alpha\beta,\alpha\gamma,\beta\gamma,\gamma^2)$ is smooth
and it is smoothable by
Corollary~\ref{ref:unionoflineandpointissmoothable:cor}.
Using Lemma \ref{lemma: homogeneous lies in every component of fiber} again,
this smooth and smoothable point lies in each irreducible component of the fiber,
so each irreducible component of the fiber is contained in $R_3^{11}$.

Now we consider the last two cases, where one finds that the homogeneous ideals
$J_3^\perp$, $J_4^\perp$ are not smooth points (although they are monomial, hence, smoothable).
So we need to develop a more detailed description of the fibers in these cases.
In Lemma \ref{lemma: 132221 x5 y4} we show that the fiber $\pi_{\vh}^{-1}(J_3^\perp)$ is contained in $R_3^{11}$
and in Lemma \ref{lemma: 132221 x5 x3y} we do the same for $J_4$.
\qed
\end{proof}

\def\augm#1{\tilde{#1}}%
\begin{lemma}\label{lemma: 132221 x5 y4}
Let $\vh=(1,3,2,2,2,1)$ and $J = \langle x^5, y^4, z \rangle$, $I = J^\perp$.
Then $\pi_{\vh}^{-1}(I) \subset R_3^{11}$.
\end{lemma}
\begin{proof}
First we will show that the fiber $\pi_{\vh}^{-1}(I)$ is irreducible,
then we will display a smooth and smoothable point in the fiber.
To begin, $I$ is generated by
$f_1 = \alpha\beta$, $f_2 = \alpha\gamma$, $f_3 = \beta\gamma$, $f_4 = \gamma^2$, $\beta^5$, $\alpha^6$.
Let $I' \in \pi_{\vh}^{-1}(I)$. Then
\begin{equation}\label{eq:idealform}
I' = (F_1, F_2, F_3, F_4, \beta^5) + (\alpha,\beta,\gamma)^{6},
\end{equation}
where $F_i = f_i + g_i$ and
each $g_i$ involves monomials of degree $3$ or greater that are not in $I$.
Those monomials are $\alpha^3, \beta^3, \alpha^4, \beta^4, \alpha^5$.
We can write, for each $i=1,2,3,4$,
\begin{equation*}
  g_i = a_i \alpha^3 + b_i \beta^3 + c_i \alpha^4 + d_i \beta^4 + e_i \alpha^5 .
\end{equation*}
This embeds the fiber $\pi_{\vh}^{-1}(I)$ into $\A^{20}$ with coordinates
$a_1,\dotsc,e_4$. It remains to find its equations, that
is, determine which ideals $I'$ of the form~\eqref{eq:idealform} have initial ideal $I$.
We claim that $\pi_{\vh}^{-1}(I)$ is defined by the equations
\begin{equation}\label{eqn: 132221 x5 y4 fiber equation}
  b_2 = a_3 = a_4 = b_4 = a_1 a_2 + c_3 = a_2^2 + c_4 = 0.
\end{equation}

Since $\init(I') \supset I$ and $\dim T/\init(I') = \dim T/I'$ we
have $\init(I') = I$ if and only if $\dim T/I = \dim T/I'$.
Consider the elements $\augm{g_i} =  a_1\alpha^3 + b_i \beta^3 + t\cdot (c_i\alpha^4 +
d_i\beta^4) + t^2\cdot e_i\alpha^5\in T[t]$ and $\augm{F_i} = f_i + t\augm{g_i}$.
Define the ideal
\begin{equation}
    \augm{I'} = (\augm{F_1}, \augm{F_2}, \augm{F_3}, \augm{F_4}, \beta^5) + (\alpha, \beta, \gamma)^{6}.
\end{equation}
Clearly, the fiber of $\augm{I'}$ over $t = 1$ is $I'$ and over $t = 0$ is $I$.
Also the family is flat over $\kk[t^{\pm 1}]$ because of the torus action.
Therefore, $\augm{I'}$ is flat
if and only if all fibers have the same length, if and only if $\dim T/I' = \dim T/I$.
That is, $I'\in \pi_{\vh}^{-1}(I)$ if and only if $\augm{I'}$ is flat.
Flatness of $\augm{I'}$ is equivalent to the following condition
(see, for example, \cite[p.~11]{artin_deform_of_sings} or \cite[Corollary 7.4.7]{MR2363237}).
\begin{equation*}
    \mbox{Every relation }\sum f_i r_i = 0\mbox{ with
    }r_i\in T\mbox{ lifts to }\sum \augm{F_i} R_i = 0\mbox{ with }R_i\in T[t].
\end{equation*}
That is, there exist $R'_i \in T[t]$ such that
$0 = \sum \augm{F_i} (r_i + t R'_i) = \sum r_i f_i + t \sum r_i \augm{g_i} + t \sum R'_i \augm{F_i}$,
equivalently $\sum r_i \augm{g_i} = - \sum R'_i \augm{F_i} \in \augm{I'}$.
So $\augm{I'}$ is flat if and only if the following holds.
\begin{equation}
    \label{it:syzygies} \mbox{For every relation }\sum f_i r_i = 0\mbox{ with
    }r_i\in T\mbox{ we have }\sum \augm{g_i} r_i \in \augm{I'}.
\end{equation}
The relations between the $f_i$ are the syzygies of $I$. They are generated
by four linear syzygies, two quartic syzygies, and two quintic syzygies (direct check). It is enough
to check~\eqref{it:syzygies} for those generators.
Since $\augm{I'} \supset (\alpha, \beta, \gamma)^6$, the
property~\eqref{it:syzygies} is automatically satisfied for quartic and quintic syzygies.
The linear generators are given by
\begin{equation*}
     \gamma f_1 = \beta f_2,
     \quad
     \beta f_2 = \alpha f_3,
     \quad
     \gamma f_2 = \alpha f_4,
     \quad
     \gamma  f_3 = \beta f_4.
\end{equation*}
By \eqref{it:syzygies}, the fiber is cut out by the conditions
\begin{equation*}
    \gamma \augm{g_1} - \beta \augm{g_2}\in\augm{I'},
    \quad
    \beta \augm{g_2} - \alpha \augm{g_3}\in\augm{I'},
    \quad
    \gamma \augm{g_2} - \alpha \augm{g_4}\in\augm{I'},
    \quad
    \gamma  \augm{g_3} - \beta \augm{g_4}\in\augm{I'}.
\end{equation*}
We now check that they unfold into \eqref{eqn: 132221 x5 y4 fiber equation}.
Consider an ideal $I'\in \pi^{-1}_{\vh}(I)$.
The element $\gamma \augm{g_1} - \beta \augm{g_2}$ lies in
$\augm{I'}$ by \eqref{it:syzygies}. Since
\begin{multline*}
    \gamma \augm{g_1} - \beta \augm{g_2}
    = a_1 \alpha^3\gamma + b_1 \beta^3\gamma - a_2 \alpha^3\beta - b_2 \beta^4 \\
    + t(c_1 \alpha^4\gamma + d_1 \beta^4\gamma - c_2 \alpha^4\beta - d_2 \beta^5)
    + t^2(e_1 \alpha^5\gamma - e_2 \alpha^5\beta)
\end{multline*}
lies in $\augm{I'}$, its initial form lies in $I$, which implies $b_2 = 0$.
Similarly, by considering the initial forms of $\gamma \augm{g_1} - \alpha
\augm{g_3} \in I'$ we deduce that $a_3 = 0$;
from $\gamma \augm{g_2} - \alpha \augm{g_4} \in I'$ we get $a_4 = 0$;
from $\gamma \augm{g_3} - \beta \augm{g_4} \in I'$ we get $b_4 = 0$.
Note the following relations:
\begin{equation*}
    \alpha^3 \beta \equiv - a_1 t\alpha^5 , \quad
    \alpha^3 \gamma \equiv -a_2 t\alpha^5 , \quad
    \beta^3 \gamma \equiv - b_3 t\beta^5 \pmod{\augm{I'}}.
\end{equation*}
Using these relations, together with $b_2 = a_3 = a_4 = b_4 = 0$, we check that
\begin{equation*}
    \beta \augm{g_2} - \alpha \augm{g_3}\equiv -t(a_1a_2 + c_3)\alpha^5
    \pmod{\augm{I'}}.
\end{equation*}
This implies that $-t(a_1a_2 + c_3)\alpha^5\in \augm{I'}$, so by evaluating at
$t =1$ we get $(a_1a_2 + c_3)\alpha^5\in I'$. Hence, the leading form
$(a_1a_2 + c_3)\alpha^5$ is in $I$. Therefore, $a_1a_2 + c_3 = 0$.
Similarly, $\gamma \augm{g_2} - \alpha \augm{g_4} \equiv -(a_2^2 + c_4) \alpha^5 \pmod{\augm{I'}}$
which gives the condition $a_2^2 + c_4 = 0$,
whereas for $\gamma \augm{g_3} - \beta \augm{g_4}$ and $\gamma \augm{g_1} - \beta \augm{g_2}$ we get trivially zero.
Thus, \eqref{eqn: 132221 x5 y4 fiber equation} is satisfied for every
$I'$ in the fiber. Conversely, the above reasoning implies that each $I'$
satisfying~\eqref{eqn: 132221 x5 y4 fiber equation} lies in the fiber.
This shows that the fiber is irreducible, in fact isomorphic to $\A^{14}$
via projection to the coordinates $a_1$, $a_2$, $b_1$, $b_3$, $c_1$, $c_2$, $d_1,\dotsc,e_4$.

Finally, let
$I' = (\alpha^6, \allowbreak \beta^5, \allowbreak \alpha\beta, \allowbreak
\alpha\gamma, \allowbreak \beta\gamma+\alpha^5, \gamma^2)$.
It is smoothable by
Corollary~\ref{ref:unionoflineandpointissmoothable:cor}.
We verify computationally that $I'$ is a smooth point.
\qed
\end{proof}

\begin{lemma}\label{lemma: 132221 x5 x3y}
Let $\vh=(1,3,2,2,2,1)$ and $J = \langle x^5, x^3y, z \rangle$, $I = J^\perp$.
Then $\pi_{\vh}^{-1}(I) \subset R_3^{11}$.
\end{lemma}
\begin{proof}
The proof directly follows the argument of Lemma~\ref{lemma: 132221 x5
y4}.
The ideal $I$ is generated by $f_1 =
\alpha\gamma$,
$f_2 = \beta^2$, $f_3 = \beta\gamma$, $f_4 = \gamma^2$,
$\alpha^4\beta$, $\alpha^6$.
Let $I' \in \pi_{\vh}^{-1}(I)$. Then
\begin{equation}
I' = (F_1, F_2, F_3, F_4, \beta^5) + (\alpha,\beta,\gamma)^{6},
\end{equation}
where $F_i = f_i + g_i$ and $g_i = a_i \alpha^3 + b_i \alpha^2 \beta + c_i
\alpha^4 + d_i \alpha^3 \beta + e_i \alpha^5$.
The syzygies among $f_i$'s are again generated by linear, quartic, and quintic syzygies. The linear generators are
$\beta f_1 - \alpha f_3$, $\gamma f_1 - \alpha f_4$, $\gamma f_2 - \beta f_3$,
$\gamma f_3 - \beta f_4$. An analysis of the resulting conditions gives the following equations for
$\pi_{\vh}^{-1}(I)$:
\begin{equation}\label{eqn: 132221 x5 x3y fiber equation}
  a_1-b_3=a_3=a_4=b_4 = a_2 b_1 + c_3 = a_1^2 + c_4 = 0.
\end{equation}
This shows that the fiber $\pi_{\vh}^{-1}(I)$ is irreducible, in fact isomorphic to $\A^{14}$
via projection to the coordinates $a_1,a_2,b_1,b_2,c_1,c_2,d_1,\dotsc,e_4$.
A smooth and smoothable point in the fiber is
\(
  I' = (\alpha^6, \alpha^4\beta, \alpha\gamma, \allowbreak \beta^2, \allowbreak \beta\gamma, \allowbreak \gamma^2+\alpha^5) .
\)
It is smoothable by Corollary~\ref{ref:unionoflineandpointissmoothable:cor}
and is computationally verified to be a smooth point.
\qed
\end{proof}

\subsection[Case (1,3,3,2,2)]{Case $\vh=(1,3,3,2,2)$}\label{sec:13322}

\begin{lemma}
Let $f \in \kk[x,y]$ be a homogeneous form of degree $d \geq 3$.
Either $f = \ell^d + m^d$ or $f = \ell^{d-1} m$ for some linear forms $\ell$ and $m$,
or else $f$ is determined up to scalar multiple by the subspace $\langle f \rangle_{d-1}$,
in the sense that if $f \neq \ell^d + m^d, \ell^{d-1}m$
and $g \in \kk[x,y]$ is a homogeneous form of degree $d$ such that $\langle f \rangle_{d-1} = \langle g \rangle_{d-1}$,
then $g$ is a scalar multiple of $f$.
\end{lemma}
\begin{proof}
If $f = \ell^d$ for a linear form $\ell$ then
$\langle g \rangle_{d-1} = \langle f \rangle_{d-1} = \Span\{\ell^{d-1}\}$, so $g$ is a scalar multiple of $\ell^d$.
Otherwise let $I = f^\perp$ and $J = g^\perp$.
The assumption $\langle f \rangle_{d-1} = \langle g \rangle_{d-1}$ means $I_{d-1} = J_{d-1}$.
Assuming $f \neq \ell^d,\ell^d+m^d,\ell^{d-1}m$ means
that $f^\perp$ has no generators of degree $\geq d$
by \cite[Proposition 1.6, Theorem 1.7]{MR3426613}.
So $I_d$ is determined by $I_{d-1} = (\langle f \rangle_{d-1})^\perp$.
Since $J_{d-1} = I_{d-1}$, these generate the same degree $d$ part, $J_d = I_d$.
But $J_d$ is perpendicular to $g$ while $I_d$ is perpendicular to $f$, so $g$ and $f$ are linearly dependent.
\qed
\end{proof}

\begin{proposition} \label{case7}
Let $\vh$ be $(1,3,3,2,2)$. Then $ H_3^{\vh} \subset R_3^{11}$.
\end{proposition}
\begin{proof}
Let $J$ be a graded inverse system in $S = \kk[x,y,z]$ with Hilbert function $(1,3,3,2,2)$.
Then $\dim J_4 = 2$, say $J_4 = \Span\{f,g\}$.
Each $f,g$ has first derivatives in $J_3$, so each $f,g$ involves at most two variables.
If both $\langle f \rangle$ and $\langle g \rangle$ have Hilbert function
$(1,*,*,1,1)$,
then $f = \ell^4$, $g = m^4$ for independent linear forms $\ell,m$, and we change coordinates
so $(f,g) = (x^4,y^4)$.
Otherwise at least one, say $\langle f \rangle$, has Hilbert function $(1,*,*,2,1)$.
Then by Proposition \ref{proposition: *k1 *11 *221} there is a coordinate change so that $f \in \kk[x,y]$.
We have $\langle g \rangle_3 \subseteq J_3 = \langle f \rangle_3 \subset \kk[x,y]$.
This shows that $\alpha \hook g, \beta \hook g, \gamma \hook g$ have no terms involving $z$.
This implies $g$ has no terms involving $z$.
So $g \in \kk[x,y]$ as well.

Now there are various cases, according as $f = \ell^4+m^4$
(which we may take to be $x^4+y^4$ after a change of coordinates),
$f = \ell^3 m$ (equivalently, $x^3 y$), or something else;
and $\dim \langle g \rangle_3 = 1$ or $2$.
In every case one checks that either
$\Span\{f,g\} = \Span\{x^4,y^4\}$ or $\Span\{f,g\} = \Span\{x^4,x^3 y\}$,
after a change of coordinates.

In either case, $J$ is generated by $J_4$, some quadratic form $Q$, and possibly linear forms:
$J$ is generated, possibly redundantly, either by $\{x^4,y^4,Q,x,y,z\}$ or by $\{x^4,x^3y,Q,x,y,z\}$,
where $Q$ is linearly independent from $\{x^2,y^2\}$ in the first case or $\{x^2,xy\}$ in the second case.

Now we claim that there is an automorphism of $S_1 = \Span\{x,y,z\}$ that takes $J$ to
one of the following.
If $J_4$ is generated by $x^4, y^4$ then we claim there is an automorphism taking $J$
to the inverse system generated by $\{x^4,y^4,Q,x,y,z\}$ where $Q \in \{z^2, z^2+xy, z(x+y),zx,xy\}$.
And if $J_4$ is generated by $x^4, x^3 y$ then we claim there is an automorphism taking $J$
to the inverse system generated by $\{x^4,x^3 y,Q,x,y,z\}$ where $Q \in \{z^2, z^2+y^2, yz, y^2+xz, y^2, xz\}$.

First suppose $J$ is generated by $x^4,y^4,Q,x,y,z$.
Write $Q = axy + bxz + cyz + dz^2$, where we can eliminate $x^2,y^2$ terms since $x^2,y^2 \in J_2$.
If $d \neq 0$ then replacing $z$ with a suitable linear combination of $z,x,y$ allows us to eliminate the $xz,yz$
terms by completing the square, as well as simultaneously rescaling $z$ to get rid of the coefficient $d$.
Then $Q = a'xy + z^2$.
If $a'=0$ then $Q = z^2$, and if $a' \neq 0$ then rescaling $x,y$ gives $Q = z^2 + xy$.
On the other hand, if $d = 0$, then rescaling $x,y,z$ allows us to get rid of the coefficients $a,b,c$, so we may assume each of them
is $0$ or $1$.
This shows $Q \in \{z^2,z^2+xy,xy+xz+yz,xy+xz,xy+yz,xz+yz,xy,xz,yz\}$.
By symmetry, interchanging $x$ and $y$ allows us to eliminate the cases $xy+yz,yz$ since these are respectively isomorphic to $xy+xz,xz$.
And replacing $z$ with $z-y$ takes $xy+xz=x(y+z)$ to $xz$.
Similarly, replacing $z$ with $z-y$ takes $xy+xz+yz$ to $xz+yz-y^2$, and $\Span\{x^2,y^2,xz+yz-y^2\} = \Span\{x^2,y^2,xz+yz\}$,
so this case is also equivalent to $Q=xz$.
This finishes the analysis of the case $J_4 = \Span\{x^4,y^4\}$.

The case $J_4 = \Span\{x^4,x^3y\}$ is similar.
Instead of a symmetry interchanging $x$ and $y$, we can replace $y$ with $y+ax$, since $\Span\{x^4,x^3y\} = \Span\{x^4,x^3(y+ax)\}$.
Write $Q = axz + by^2 + cyz + dz^2$, after eliminating $x^2, xy$ terms.
If $d\neq0$ then a substitution for $z$ eliminates $xz,yz$ terms, yielding $Q = b'y^2 + z^2$.
Rescaling $y$ if necessary, $Q = z^2$ or $Q = y^2 + z^2$.
If $d = 0$ then rescaling $x,y,z$ to eliminate the $a,b,c$ coefficients gives
$Q \in \{xz + y^2 + yz, xz+y^2, xz+yz, y^2+yz, xz, y^2, yz\}$.
Appropriate substitutions for $y$ and $z$ take the cases $xz+y^2+yz, xz+yz, y^2+yz$ all to $yz$.

Now by Lemma~\ref{irreducibleFiber} each fiber over a point in $\cH_3^{\vh}$ is irreducible.
Thus, it suffices to find a smooth and smoothable inverse system $J'$ such that $\lead(J') = J$ for each of the normal forms $J$.
For the case that $J_4$ is spanned by $x^4$ and $y^4$ see Table \ref{table: 13322 deformations x4 y4}.
For the case that $J_4$ is spanned by $x^4$ and $x^3 y$ see Table \ref{table: 13322 deformations x4 x3y}.
\qed
\end{proof}

\begin{table}
\begin{equation*}
\begin{array}{lll}
  \toprule
  Q & J' & \text{deformation of ideal of $J'$} \\
  \midrule
  z^2 + xy & \langle x^4, y^4, z^2+xy \rangle            & (\beta\gamma,2\alpha\beta-\gamma^2,\alpha\gamma,\alpha^5,\beta^5+t\beta^4) \\
  z^2      & \langle x^4+x^2y+x^2z+z^3,y^4,z^2 \rangle   & (\beta\gamma,\alpha\beta-\alpha\gamma,3\alpha^2\gamma-\gamma^3,\alpha^3-12\alpha\gamma,\beta^5+t\beta^4) \\
  xy       & \langle x^4+x^2y+x^2z+xy^2,y^4,xy,z \rangle & (\gamma^2,\beta\gamma,\alpha\beta^2-\alpha^2\gamma,\alpha^2\beta-\alpha^2\gamma,\alpha^3-12\alpha\gamma,\beta^5+t\beta^4) \\
  xz       & \langle x^4+xz^2,y^4,xz \rangle             & (\beta\gamma,\alpha\beta,\alpha^2\gamma,\alpha^3-12\gamma^2,\beta^5+t\beta^4) \\
  (x+y)z   & \langle x^4,y^4,xz+yz \rangle               & (\gamma^2,\alpha\gamma-\beta\gamma,\alpha\beta,\beta^5,\alpha^5+t\alpha^4) \\
  \bottomrule
\end{array}
\end{equation*}
\caption{Smooth and smoothable inverse systems $J'$ with Hilbert function $(1,3,3,2,2)$
and $\lead(J')_4$ spanned by $x^4$, $y^4$.} \label{table: 13322 deformations x4 y4}
\end{table}

\begin{table}
\begin{equation*}
\begin{array}{lll}
  \toprule
  Q & J' & \text{deformation of ideal of $J'$} \\
  \midrule
  z^2     & \langle x^4,x^3y+x^2z+z^3,z^2 \rangle        & (\beta\gamma,\beta^2,3\alpha^2\gamma-\gamma^3,\alpha^2\beta-3\alpha\gamma,\alpha^5+t\alpha^4) \\
  z^2+y^2 & \langle x^4,x^3y,y^2+z^2 \rangle             & (\beta\gamma,\alpha\gamma,\beta^2-\gamma^2,\alpha^4\beta,\alpha^5+t\alpha^4) \\
  yz      & \langle x^4,x^3y+x^2z,yz \rangle             & (\gamma^2,\beta^2,\alpha\beta\gamma,\alpha^2\beta-3\alpha\gamma,\alpha^3\gamma,\alpha^5+t\alpha^4) \\
  y^2+xz  & \langle x^4+2x^2z,x^3y+xyz,y^2+xz \rangle     & (\gamma^2+t\gamma,\beta^2\gamma,\beta^3,\alpha\beta^2,\alpha^2\beta-6\beta\gamma,\alpha^3-6\alpha\gamma+3\beta^2) \\
  y^2     & \langle x^4+2x^2z,x^3y+x^2z+xyz,y^2,z \rangle & (\gamma^2,\beta^2\gamma,\beta^3+t\beta^2,\alpha\beta^2,\alpha^2\beta-6\beta\gamma,\alpha^3-6\alpha\gamma+12\beta\gamma) \\
  xz      & \langle x^4,x^3y+x^2z+xz^2,xz \rangle        & (\beta\gamma,\beta^2,\alpha^2\gamma-\alpha\gamma^2,\alpha^2\beta-3\gamma^2,\alpha^5+t\alpha^4) \\
  \bottomrule
\end{array}
\end{equation*}
\caption{Smooth and smoothable inverse systems $J'$ with Hilbert function $(1,3,3,2,2)$
and $\lead(J')_4$ spanned by $x^4$, $x^3 y$.}\label{table: 13322 deformations x4 x3y}
\end{table}

\subsection[Case (1,3,3,3,1)]{Case $\vh = (1,3,3,3,1)$}\label{sec:13331}

First we consider separately a special case, where the quadrics in the inverse systems
have a most special form.
\begin{proposition}\label{ref:13331quadricsinxy:prop}
    Let $\vh = (1, 3, 3, 3, 1)$, let $I \in H_3^{\vh}$ and let $J$ be its
    inverse system. Suppose $x^2, xy, y^2 \in J$. Then $I$ is contained in the
    smoothable component.
\end{proposition}
\begin{proof}
    The inverse system $J$ has a quartic generator and its leading form $f$ is
    uniquely determined. Since $\vh(2) = 3$ and $x^2, xy, y^2\in \lead(J)$, we
    see that $f\in \kk[x, y]$. We consider two cases. In each case we show
    that the space of possible $J$ is irreducible and find a smooth and
    smoothable point there.

            First suppose $f$ is annihilated by a linear form in $\kk[\alpha, \beta]$.
            Then, up to coordinate change, we have $f = x^4$.
            Consider the family of tuples
            $
              (x^4 + c + q, \allowbreak c_1 + q_1, \allowbreak c_2+q_2, \allowbreak x, \allowbreak y, \allowbreak z),
            $
            where $c_i$, $c$ are cubics and $q_i$, $q$
            are quadrics, with the condition that ${\gamma\hook c},
            {\beta\hook c}$ lie in $\Span\{
            x^2, xy, y^2\}$ and also all derivatives of $c_i$ lie in
            $\Span\{x^2, xy, y^2\}$. The space of polynomial tuples satisfying
            these conditions is an affine space.
            Each inverse system $K$ generated by a tuple as above has Hilbert
            function at most $(1, 3, 3, 3, 1)$. Thus, a \emph{general} one
            has Hilbert function exactly $(1, 3, 3, 3, 1)$. Denote the
            irreducible family of such $K$'s by $\mathcal{F}$.
            Then
            $\mathcal{F}$ gives a morphism to the Hilbert scheme
            $H_3^{\vh}$ and the image contains $J$. The image contains also
            $J_0 = \langle x^4+x^2z , x^2y , xy^2 ,x,y,z \rangle$.
            A deformation of its ideal is given by
            \(
              (\beta\gamma,\gamma^2+t\gamma,\beta^3,\alpha^3-12\alpha\gamma,\alpha^2\beta^2).
            \)
            For $t \neq 0$ this is supported at more than one point, hence, $J_0^\perp$ is
            smoothable. And $J_0^\perp$ is smooth as well, hence, the whole image of
            $\mathcal{F}$ is contained in $R_3^{11}$ by Lemma~\ref{lemma:
            smooth and smoothable is enough}.

Suppose now $f$ is not annihilated by a linear form in
            $\kk[\alpha, \beta]$. Then the proof of the previous case applies
            with the difference we consider the family of $g + c + q$, $c_1 + q_1$
            where $g \in \kk[x,y]_4$
            with the condition that $\gamma\hook c$ and all derivatives of $c_1$ lie in
            $\Span\{x^2, xy, y^2\}$. The smooth and smoothable point is
            given by the inverse system
            $\langle x^2y^2+xyz,x^3,z \rangle$ and a deformation of the corresponding ideal is given by 
            $(\gamma^2,\beta^2\gamma,\alpha^2\gamma,\beta^3,\alpha\beta^2-4\beta\gamma,\alpha^2\beta-4\alpha\gamma,\alpha^4+\alpha^3t)$.
\qed
\end{proof}

\begin{proposition}\label{ref:13331quadricsgeneral:prop}
Let $\vh = (1,3,3,3,1)$ and let J be a graded inverse system with Hilbert function $\vh$.
    Then up to coordinate change $J_2$ is the span of one of the following sets:
    \begin{equation*}
        \{x^2, xy, y^2\}, \quad \{ x^2, y^2, z^2\}, \quad \{x^2,yz, z^2\},
        \quad \{xz,yz,z^2\},\quad \{x^2+yz,xz,z^2\}.
    \end{equation*}
 Let $A = T/(J_2^\perp)$. Then $\dim A_3 = 4$ if $J_2 = \Span\{x^2,xy,y^2\}$ and $\dim A_3 = 3$ otherwise.
\end{proposition}

\begin{proof}
Let $I = (J_2^\perp)$.
By Macaulay's bound, $\dim A_3 \leq (\dim A_2)^{\langle 2 \rangle} = 3^{\langle 2 \rangle} = 4$.
If $\dim A_3 = 4$, then by Lemma~\ref{Quadrics via Gotzmann}
the quadrics in $I_2$ share a common linear factor.
After a change of coordinates, $I_2$ is spanned by $\{\gamma^2, \beta\gamma,
\alpha\gamma\}$.
Then $J$ contains the quadrics $x^2, xy, y^2$.

So we reduce to the case $\dim A_3 = 3$, equivalently $\dim I_3 = 7$.
We start with the claim that,
when $I = (I_2)$ is generated by $3$ quadrics and $\dim I_3 = 7$,
then \emph{$I$ is the saturated ideal of a zero-dimensional degree $3$ scheme
in $\PP^2$}.
The space $T_1 \otimes I_2$ has dimension $9$ and maps by multiplication surjectively to $T_1 \otimes I_2 \to I_3$,
so the kernel has dimension $2$, which means there are $2$ linear syzygies among the quadrics in $I_2$.
The minimal free resolution of $A = T/I$ is equal to
\begin{equation*} 
  0 \leftarrow T \leftarrow T(-2)^{\oplus 3} \leftarrow T(-3)^{\oplus 2}
  \oplus T(-4)^{\oplus q} \oplus F'
    \leftarrow T(-4)^{\oplus p} \oplus F'' \leftarrow 0,
\end{equation*}
where $F', F''$ are sums of $T(-i)$ with $i > 4$.
We will show that $p = q = 0$.
First, if $p = \beta_{3,4}(I) \neq 0$ then
$I$ contains the ideal $\ell(\alpha,\beta,\gamma)$ for some linear form $\ell$,
by \cite[discussion following Theorem 8.15, p.~162]{eisenbud:syzygies}.
But this is the case $\dim A_3 = 4$ which we have already treated.
Since we are now assuming $\dim A_3 = 3$, then we must have $p=0$.
Next, we compute $\dim A_4$ by considering the free resolution above:
\begin{equation*}
  \dim A_4 = \dim T_4 - 3 \dim T_2 + (2 \dim T_1 + q \dim T_0) - 0,
\end{equation*}
where the final $0$ reflects $p=0$.
This gives
$\dim A_4 = 15 - 3 \cdot 6 + 2 \cdot 3 + q = 3+q$.
At the same time, $3+q = \dim A_4 \leq (\dim A_3)^{\langle 3 \rangle} = 3^{\langle 3 \rangle} = 3$.
So $q = 0$.

Now $I$ is generated in degree $2$ and $\dim A_4 = (\dim A_3)^{\langle 3 \rangle} = 3$.
By Gotzmann's Persistence Theorem, $\dim A_k = 3$ for all $k \geq 3$.
This shows $Z = \Proj A$ has Hilbert polynomial $3$, so $Z = V(I) \subset \PP^2$ is zero-dimensional and has degree $3$.
To see that $I$ is saturated, let $I'$ be the saturation of $I$.
Since the quadrics in $I$ share no common linear factor, $Z$ is not contained in any line,
so $I'$ contains no linear forms.
Then $\dim (T/I')_1 = 3$.
The Hilbert function of a saturated ideal is nondecreasing, so for every $k \geq 1$,
\(
  3 = \dim (T/I')_1 \leq \dim (T/I')_k \leq \dim (T/I)_k = 3,
\)
which shows $I' = I$.
This completes the proof of the claim that $I$ is the saturated ideal of a degree $3$ zero-dimensional scheme $Z$ in $\PP^2$.

Since $Z$ is cut out by quadrics, the intersection of $Z$ with any line has degree at most $2$.
Then $Z$ is one of the following.
\begin{enumerate}
\item $Z$ may be a disjoint union of three non-collinear reduced points.
We change coordinates so that $Z = \{[1:0:0], [0:1:0], [0:0:1]\}$. Then $I_2 = \Span\{\alpha\beta,\alpha\gamma,\beta\gamma\}$.
\item $Z$ may be the union of a reduced point with a zero-dimensional scheme of degree $2$.
We choose coordinates so that the reduced point is $[1:0:0]$ and the scheme of degree $2$
is supported at $[0:0:1]$ and is contained in the line spanned by $[0:0:1]$ and $[0:1:0]$.
Then $I_2 = \Span\{\alpha\beta,\alpha\gamma,\beta^2\}$.
\item $Z$ may be a scheme of degree $3$ supported at a point which we may take to be $[0:0:1]$. Then after a change of coordinates either 
 $I_2 = \Span\{\alpha^2,\alpha\beta,\beta^2\}$ or $I_2 = \Span\{\alpha^2 - \beta\gamma, \alpha\beta, \beta^2\}$.
\end{enumerate}
To see the last claim, first note that if $q = \ell_1 \ell_2 \in I_2$ is a reducible quadric
then both components $\ell_1,\ell_2$ pass through $[0:0:1]$, because $Z$ is not contained in any single line.
If every quadric in $I_2$ is reducible then $I_2 = \Span\{\alpha^2,\alpha\beta,\beta^2\}$ consists of all the quadrics
that are singular at $[0:0:1]$.
Otherwise, there is a smooth quadric in $I_2$ and we can choose coordinates so that it is $q = \alpha^2-\beta\gamma$.
A quadric $q' \in I_2$ intersects $q$ in $Z$ plus one more point.
If the extra point is also $[0:0:1]$ then $q' = \beta^2 + \lambda q$ for some scalar $\lambda$.
Otherwise $q' = \ell\beta + \lambda q$ where $\ell$ is the line through $[0:0:1]$ and the extra point.
So $I_2$ is spanned by $q$, $\alpha\beta$, and $\beta^2$.
Now, having the normal forms of $Z$, we calculate $J_2$ as
$I_2^{\perp}$, obtaining the list above.
\qed
\end{proof}

\begin{lemma}\label{ref:quadricsandleadingisenough:lem}
    Fix a three dimensional space of quadrics $Q$ and a subspace $A \subset
    \langle \alpha, \beta, \gamma\rangle$.
    Suppose that the derivatives of $Q$ span $\langle x, y, z\rangle$.
    Consider the set $\mathcal{J}(Q, A)$ of inverse
    systems $J$ such that
    \begin{enumerate}
        \item $J$ has Hilbert function $(1, 3, 3, 3, 1)$,
        \item $Q$ equals $J_2$,
        \item $A$ is equal to the space of linear forms annihilating the
            quartic in $\lead(J)$.
    \end{enumerate}
    \label{ref:13331irreducibilityLemma}
    Then $\mathcal{J}(Q, A)$ is irreducible or empty.
\end{lemma}
\begin{proof}
    Suppose $\mathcal{J}(Q, A)$ is non-empty. Let $I = Q^{\perp}$.
    Let $\vh$ be the Hilbert function of $T/I$.
    By Proposition~\ref{ref:13331quadricsgeneral:prop}, we have either $Q =
    \Span\{x^2, xy, y^2\}$ up to coordinate change or $\vh(3) = 3$. The
    equality $Q = \Span\{x^2, xy, y^2\}$ is impossible, since derivatives of
    $Q$ span $x, y, z$. Thus, $\vh(3) = 3$.
    The remaining part of the proof resembles the proof of
    Proposition~\ref{ref:13331quadricsinxy:prop}.
    Let $a = \dim A$.
    Consider the set $\mathcal{F}$ of tuples of polynomials
    \begin{equation}
        f + c + q,\ c_1 + q_1, \ldots , c_a + q_a\in S,
    \end{equation}
    such that
    \begin{enumerate}
        \item $f$ is homogeneous of degree four and annihilated by $A$ (and
            possibly other linear forms),
        \item all $c_i$ and $c$ are homogeneous of degree three,
        \item all $q_i$ and $q$ are homogeneous of degree two,
        \item both $f$ and all $c_i$ are annihilated by $I$,
        \item the space $A\hook c$ is contained in $Q$.
    \end{enumerate}
    All given conditions are linear in coefficients of polynomials, thus,
    $\mathcal{F}$ is an affine space.
    Consider an open (possibly empty) subset $\mathcal{F}_0 \subset
    \mathcal{F}$ consisting of tuples where $f$ is annihilated
    exactly by $A$ and such that $c_i$ are linearly independent and
    $\Span\{c_i\}$ is disjoint from the space of partial derivatives of
    $f$. Then $\mathcal{F}_0$ is irreducible as an open set in affine space.

    Consider an inverse system $J$ generated by all linear forms and a tuple in $\mathcal{F}_0$. We
    now prove that its Hilbert function $\tilde{\vh}$ is at most $(1, 3, 3, 3, 1)$
    position-wise. By Proposition~\ref{proposition: inverse system same
    hilbert function as leading forms}, it is enough to show that the Hilbert
    function of $\lead(J)$ is
    at most $(1, 3, 3, 3, 1)$.
    It is clear that $\tilde{\vh}(0)=\tilde{\vh}(4)=1$ and $\tilde{\vh}(1)\leq 3$.
    All cubic terms in $\lead(J)$ are leading forms of combinations of $c_i$ and partials of $f$. Thus, they are
    annihilated by $I$. The space of cubics annihilated by $I$ is $\vh(3) =
    3$-dimensional, thus, $\tilde{\vh}(3) \leq 3$.
    Consider now the quadrics in $\lead(J)$. They are combinations of leading
    forms of partials of $f$, of $c_i$ and also of $A\hook c$. All those forms
    lie in $Q$, thus, $\tilde{\vh}(2)\leq 3$. Therefore, $\tilde{\vh} \leq (1,
    3, 3, 3, 1)$ position-wise.

    Since $(1, 3, 3, 3, 1)$ is the maximal possible value of $\tilde{\vh}$, the set
    $\mathcal{F}_{gen} \subset \mathcal{F}_0$ consisting of systems with
    Hilbert function $(1, 3, 3, 3, 1)$ is open, thus,
    irreducible. It gives a map to the Hilbert scheme whose image is
    $\mathcal{J}(Q, A)$, which is, therefore, irreducible as well.
\qed
\end{proof}

\begin{proposition}
Let $\vh$ be $(1,3,3,3,1)$. Then $H_3^{\vh} \subset R_3^{11}$.
\end{proposition}
\begin{proof}
Let $J$ be a graded inverse system in $S = \kk[x,y,z]$ with Hilbert function $\vh = (1,3,3,3,1)$.
It has a unique degree $4$ generator $f$.
We will subdivide the cases according to the Hilbert function of the inverse system $K$ generated by $f$.
The Hilbert function is symmetric.
Using Macaulay's bound, we find that there are four different cases for the Hilbert function of $K$:
$(1,1,1,1,1)$, $(1,2,2,2,1)$, $(1,2,3,2,1)$, $(1,3,3,3,1)$.

Case $(1,3,3,3,1)$:
$f$ generates the entire module, hence, $J$ is Gorenstein.
And every $J'$ in the fiber $\pi_{\vh}^{-1}(J)$ is also Gorenstein.
By~\cite[Proposition~2.2]{kleppe_roig_codimensionthreeGorenstein}
or~\cite[Corollary 4.3]{kleppe_pfaffians}, any Gorenstein subscheme of $A^3$ is smoothable,
hence, all such points $J'$ lie in the smoothable component of $H_3^{11}$.

Case $(1,2,3,2,1)$:
Since $f$ has two independent first derivatives, we know it depends on only two variables, say $x,y$.
Since it spans a three-dimensional set of second derivatives, this will have
to be $\langle x^2,xy,y^2 \rangle$, so we are done by
Proposition~\ref{ref:13331quadricsinxy:prop}.

Case $(1,2,2,2,1)$:
Let $Q$ be the space of quadrics inside $J$.
By Proposition~\ref{ref:13331quadricsinxy:prop}, we may assume that
$Q\neq\Span\{x^2, xy, y^2\}$ up to coordinate change, so that derivatives of
$Q$ span $x, y, z$. Then by Proposition~\ref{ref:13331quadricsgeneral:prop}
and Lemma~\ref{ref:quadricsandleadingisenough:lem} we see
that the irreducible strata are determined by $Q$ equal to
\begin{equation}\label{eq:Qcases}
    \quad\Span\{ x^2, y^2, z^2\},\quad \Span\{x^2, yz, z^2\},\quad
    \Span\{xz,yz,z^2\},\quad \Span\{x^2+yz,xz,z^2\}
\end{equation}
and the linear forms annihilating $f$ (up to simultaneous coordinate change).
It remains to check which annihilators are possible for each $Q$.
Let $A = T/Q^{\perp}$.
By the proof of
Proposition~\ref{ref:13331quadricsgeneral:prop}, this is a homogeneous coordinate ring of a
zero-dimensional subscheme of $\Proj T$. Note that if a linear form $\sigma$
annihilates $f$, then the intersection of $\Proj A$ with the projective line $(\sigma =
0)$ has degree at least two.
We directly check that for the four cases in \eqref{eq:Qcases} we
get the following possible annihilators.
\begin{enumerate}
    \item $\alpha$ or $\beta$ or $\gamma$ for $Q = \Span\{ x^2, y^2, z^2\}$,
    \item $\alpha$ or $\beta$ for $Q = \Span\{x^2, yz, z^2\}$,
    \item $\lambda_1 \alpha + \lambda_2 \beta$ with $\lambda_i\in \kk$
        arbitrary for $Q = \Span\{xz,yz,z^2\}$,
    \item $\beta$ for $Q = \Span\{x^2+yz,xz,z^2\}$.
\end{enumerate}
Note that for first, third and fourth case there is a unique choice up to coordinate
change. Therefore, we have five distinct cases in total.
The corresponding smooth and smoothable points are presented in Table
\ref{table: 13331 12221}.
\begin{table}
\begin{equation*}
\begin{array}{ll ll}
  \toprule
  Q & (f^{\perp})_1\phantom{m} & J' & \text{deformation of ideal of $J'$} \\
  \midrule
  \Span\{x^2, y^2, z^2\} & \gamma & \langle x^4+y^4,z^3 \rangle            & (\beta\gamma,\alpha\gamma,\alpha\beta,\gamma^4+t\gamma^3,\alpha^4-\beta^4) \\
  \Span\{x^2, yz, z^2\}   & \alpha & \langle yz^3, x^3 \rangle   &      \mbox{monomial ideal, hence, smoothable}  \\
  \Span\{x^2, yz, z^2\}   & \beta  & \langle x^4+z^4,yz^2 \rangle   & (\beta^2+t\beta,\alpha\gamma,\alpha\beta,\beta\gamma^3,\alpha^4-\gamma^4) \\
  \Span\{xz,yz, z^2\}     & \beta & \langle xz^3,yz^2 \rangle             &
  \mbox{monomial ideal, hence, smoothable} \\
  \Span\{ x^2+yz, xz,z^2\}  & \beta & \langle xz^3,x^2z+yz^2 \rangle & (\beta^2,\alpha\beta,\alpha^2-\beta\gamma,\beta\gamma^3,\gamma^4+t\gamma^3) \\
  \bottomrule
\end{array}
\end{equation*}
\caption{Smooth and smoothable points $J'$ with Hilbert function $(1,3,3,3,1)$
such that the inverse system generated by $f \in J'_4$ has Hilbert function
$(1,2,2,2,1)$}\label{table: 13331 12221}
\vspace{-1em}
\end{table}

Case $(1,1,1,1,1)$: the argument is completely analogous to the previous case
up to the point where we determine possible $(f^{\perp})_1$
depending on $Q$. As before, let $A = T/Q^{\perp}$.
Note that an
annihilator $(\sigma_1, \sigma_2)$ is possible if and only if $l^4\in J$,
equivalently $[l]\in \Proj A$, where $l$ is the
linear form annihilated by $(\sigma_1, \sigma_2)$ and $[l]\in \Proj T$ is its
class.
Based on this observation, the possible annihilators in the four cases
in \eqref{eq:Qcases} are
\begin{enumerate}
    \item $(\alpha, \beta)$ or $(\beta,\gamma)$ or $(\alpha, \gamma)$ for $Q = \Span\{ x^2, y^2, z^2\}$,
    \item $(\alpha, \beta)$ or $(\beta, \gamma)$ for $Q = \Span\{x^2, yz, z^2\}$,
    \item $(\alpha, \beta)$  for $Q = \Span\{xz,yz,z^2\}$,
    \item $(\alpha, \beta)$ for $Q = \Span\{x^2+yz,xz,z^2\}$.
\end{enumerate}
The three possibilities for $Q = \Span\{ x^2, y^2, z^2\}$ are equivalent, thus,
we get five cases in total.
The list of smooth
and smoothable points is presented in Table \ref{table: 13331
11111}.
This completes the proof.
\qed

\begin{table}
\begin{equation*}
\begin{array}{llll}
  \toprule
  Q & (f^{\perp})_1\phantom{m} & J & \text{deformation} \\
  \midrule
  \{x^2, y^2, z^2\}     & (\alpha, \beta)       & \langle z^4+xz+xy+yz,x^3,y^3 \rangle
    & (\alpha\gamma-\beta\gamma,\alpha\beta-\beta\gamma,\alpha^2\beta,\alpha^4+t\alpha^3,\beta^4,\gamma^4-24\alpha\beta) \\
  \{x^2, yz, z^2\}    & (\alpha, \beta)      & \langle
    z^4,x^3+y^2,yz^2 \rangle
    & (\alpha\gamma, \alpha\beta, \beta^2\gamma, \alpha^3 +t\alpha^2- 3\beta^2,
    \beta\gamma^3, \gamma^5) \\
  \{x^2, yz, z^2\}    & (\beta, \gamma)      & \langle x^4+y^2,yz^2,z^3 \rangle
    & (\alpha\gamma,\alpha\beta,\beta^2\gamma,\beta\gamma^3,\gamma^4+t\gamma^3,\alpha\gamma^3,\alpha^4-12\beta^2) \\
  \{xz,yz, z^2\}  &  (\alpha, \beta)     & \langle y^2z+z^4,yz^2,xz^2 \rangle
    & (\alpha\beta, \alpha^2-t\alpha, \gamma^3-12\beta^2, \beta^3, \beta^2\gamma^2) \\
  \{x^2+yz, xz,z^2\}  & (\alpha, \beta)   & \langle z^4+y^2,yz^2+x^2z,xz^2 \rangle
    & (\alpha\beta,\alpha^2-\beta\gamma,\beta^3,\alpha\gamma^3,\beta\gamma^3,\gamma^4+t\gamma^3-12\beta^2) \\
  \bottomrule
\end{array}
\end{equation*}
\caption{Smooth and smoothable points $J'$ with Hilbert function $(1,3,3,3,1)$
such that the inverse system generated by $f \in J'_4$ has Hilbert function
$(1,1,1,1,1)$.}\label{table: 13331 11111}
\vspace*{-1em}
\end{table}
\end{proof}
\goodbreak

\subsection[Case (1,3,3,2,1,1)]{Case $\vh = (1, 3, 3, 2, 1, 1)$}\label{sec:133211}

\begin{proposition}
Let $\vh = (1,3,3,2,1,1)$.
Then $H_3^{\vh} \subset R_3^{11}$.
\end{proposition}
\begin{proof}
Consider a local ideal $I$ with Hilbert function $\vh= (1, 3, 3, 2, 1, 1)$,
let $J$ be the inverse system of $I$,
choose generators of $J$, and let $f$ be the generator of $J$ of degree five.
Let $\vh_{f}$ be the Hilbert function of the algebra $T/f^\perp$ apolar to $f$.
The case decomposes into five subcases, depending on $\vh_{f}$.
Since $f^{\perp} \supset I$, $T/f^\perp$ is a quotient of $A = T/I$, so that
$\vh_{f} \leq (1, 3, 3, 2, 1, 1)$.
If $\vh_{f}(3)=2$ then $\vh_{f}(1),\vh_{f}(2) \geq 2$ by Corollary~\ref{corollary: hilbert function nonincreasing}.

\begin{enumerate}
    \item $\vh_{f} = (1, a, b, 1, 1, 1)$.
        Here we argue exactly as in the proof of
        Proposition~\ref{ref:tailsofones:prop}, so we omit some details below.
        After a non-linear change of coordinates as in
        Section~\ref{sect: nonlinear}, we may assume $f = x^5 + g$ with
        $\alpha^2 \hook g = 0$.
        Since $J$ is generated by $f$ together with elements of degree $3$ or less,
        $\alpha^3\beta$ and $\alpha^3\gamma$ annihilate all of $J$.
        If $q \in (\beta,\gamma)$ is such that $\alpha^c-q$ annihilates $f$,
        then $c \geq 4$; hence, $\alpha^3 \notin I+(\beta,\gamma)$.
        Now
        Corollary~\ref{ref:algebraicversion} proves that the element is
        cleavable, hence, smoothable by Lemma \ref{lemma: cleavable implies smoothable}.
    \item $\vh_{f} = (1, 3, 3, 2, 1, 1)$. In this case $I = f^{\perp}$ and
        so $A$ is Gorenstein, hence, smoothable
        by~\cite[Proposition~2.2]{kleppe_roig_codimensionthreeGorenstein}
        or~\cite[Corollary 4.3]{kleppe_pfaffians}.
    \item $\vh_{f} = (1, 2, 3, 2, 1, 1)$. After changing coordinates so that
        $f \in \kk[x,y]$ we have $J = \langle f,
        z\rangle$, so $A$ is smoothable by
        Corollary~\ref{ref:unionoflineandpointissmoothable:cor}.
    \item $\vh_{f} = (1, 3, 2, 2, 1, 1)$. In this case after a non-linear
        change of coordinates we have $f = g + z^2$ for $g\in \kk[x, y]$ with
        Hilbert function $(1, 2, 2, 2, 1, 1)$,
        see~\cite[Proposition~4.5, Example~4.6]{cjn13}. The set of those $g$
        is irreducible by a result of Iarrobino~\cite[Proposition~4.8]{cjn13}.
        Thus, the set of $f$ is also irreducible. $J$
        is generated by $f$ and a quadric $q$ which may be chosen arbitrarily,
        modulo $\langle f \rangle_2$, thus, the set of pairs $(f,q)$ is irreducible.
        It is now enough to find a smooth, smoothable point.
        Such a point is given by the ideal $I = \langle x^5 + y^4 + z^2, \allowbreak xz \rangle^{\perp}
        = \lim_{t \to 0} (\beta\gamma, \allowbreak \alpha\beta, \allowbreak \alpha\beta^2, \allowbreak
        \alpha^2\beta, \allowbreak \beta^4+t\beta^3-12\gamma^2, \allowbreak \alpha^5-60\gamma^2)$.
    \item $\vh_{f} = (1, 2, 2, 2, 1, 1)$. As in the previous case, after a nonlinear
        change of coordinates we get $f\in \kk[x, y]$ and the set of $f\in
        \kk[x, y]$ is irreducible by~\cite[Proposition~4.8]{cjn13}. $J$
        is generated by $f$ and a quadric $q$ with no relations.
        For general $q$ after adding a multiple of $q$ to $f$ we get $f_{\text{new}}$ with
        $\vh_{f_{\text{new}}} = (1, 3, 2, 2, 1, 1)$ and, thus, reduce to the previous
        case.\qed
\end{enumerate}
\end{proof}

\subsection[Case (1,3,6,1)]{Case $\vh = (1, 3, 6, 1)$}\label{sec:1361}

\begin{proposition}\label{proposition: 1361}
Let $\vh = (1,3,6,1)$.
Then $H_3^{\vh} \subset R_3^{11}$.
\end{proposition}

    This case is in contrast with previous ones. First, it is easy to check that
    $H_3^{\vh} = \cH_3^{\vh}$,
    that is, every local ideal with Hilbert function $\vh$ is homogeneous.
    And second, it is also easy to check that $\cH_3^{\vh}$
    is irreducible, in fact a $\PP^9 = \Gr(9, 10)$. The isomorphism is
    given by sending an ideal to all its cubic equations.
    However, in this set there seem to be no smooth points.
    We argue by showing that general points in $\cH_3^{\vh}$ are smoothable,
    then by irreducibility so are all points.
    
    \begin{proof}
    Let $I \in H_3^{\vh}$ with inverse system $J$.
    Necessarily $I$ and $J$ are homogeneous.
    Let $f \in J$ be the cubic generator, which is unique up to scalar.
    We see that $H_3^{\vh} = \cH_3^{\vh}$ is parametrized by the point $[f]$ in the projective space
    of cubics in three variables, a $\PP^9$, so once again $H_3^{\vh}$ is irreducible.
    
    The cubic $f$ is the equation of a plane cubic curve.
    Assume it is a smooth curve.
    Then we may change coordinates so that $f$ is in Hesse normal form,
    that is $f = x^3 + y^3 + z^3 + 6 h xyz$ for some $h$,
    see for example \cite[Section 3.1.2]{MR2964027}.
    
    Now $J = \langle f, S_2 \rangle$. We may directly compute
    \begin{equation*}
      I = J^\perp =
      (\alpha\beta^2,\alpha^2\beta,\alpha\gamma^2,\alpha^2\gamma,\beta\gamma^2,\beta^2\gamma,\alpha^3-\beta^3,\alpha^3-\gamma^3,\alpha\beta\gamma-h\gamma^3).
    \end{equation*}
    Corollary~\ref{ref:algebraicversion} with $c = 2$
    implies that $I$ is cleavable, hence, it is smoothable. Thus, all $I$
    corresponding to smooth curves are smoothable.
    By irreducibility, so are all $I \in H_3^{\vh}$.
    \qed
    \end{proof}

\section{Smoothability of very compressed algebras}\label{sec:upto95}

    In this section we prove Proposition~\ref{ref:upto95:prop}, in other words we
    show that each element of $\cH^{\max, d}$ is smoothable. We begin
    by noting that the scheme is irreducible.
    \begin{lemma}
        The locus $\cH^{\max, d}$ is irreducible.
        \label{ref:upto95irreducibility}
    \end{lemma}
    \begin{proof}
        A scheme in
        $\cH^{\max, d}$ is uniquely determined by choice of $I$ inside
        $\mathfrak{m}^s$ and containing $\mathfrak{m}^{s+1}$. Hence,
        $\cH^{\max, d}$ is isomorphic to a Grassmannian, thus, irreducible.
    \qed
    \end{proof}
    To check smoothability we verify that a general point of the stratum is obtained as
    a \emph{$\kk^*$-limit}, a notion which we now explain.
    The scaling (homothety) action of $\kk^*$ on $\A^3$
    extends to an action
    on $\PP^3$. Take a set $\Gamma$ of $d$ points in $\PP^3$.
    For every $t\in \kk^*$ we may take $t\cdot \Gamma$.
    The \emph{$\kk^*$-limit} of $\Gamma$ is
    \(
        \Gamma' = \lim_{t\to 0}(t\cdot \Gamma).
    \)
    This is a flat limit, in the sense of~\cite[Proposition III.9.8]{HarAG}. It is
    constructed as follows. Take the graph of the $\kk^*$-action, which is a
    family
    \(
        Z^\circ_{\Gamma} \subset \kk^* \times \PP^{3},
    \)
    whose fiber over $t\in \kk^*$ is $t\Gamma$. This family is just the union
    of $n$ lines in $\kk^* \times \PP^3$ through the points $(1, p)$,
    where $p\in \Gamma$. All its fibers are isomorphic to $\Gamma$ and it is
    flat over $\kk^*$.  Let $Z_{\Gamma} \subset \kk \times \PP^3$ be the closure of $Z^\circ_{\Gamma}$.
    This family is flat
    over $\kk$, see~\cite[Proposition III.9.8]{HarAG}. Finally let
    \(
        \Gamma' = Z_{\Gamma} \cap (t = 0).
    \)
    By construction, $\Gamma'$ is smoothable (as a limit of $\Gamma$) and
    $\kk^*$-invariant.

    A general set $\Gamma$ of $d$ points imposes independent conditions on forms, hence, the
    ideal defining the limit scheme has no small-degree generators. For example,
    for $d=11$ the algebra $\Gamma'$ has Hilbert function $\vh = (1, 3, 6,
    1)$.
    After restricting to general $\Gamma$'s the $\kk^*$-limit can be made
    relative~\cite[Proof of Lemma~5.4]{CEVV} and we get a rational map
    \begin{equation}\label{eq:dominating}
        \varphi_{d}:R^d \PP^3 \dashrightarrow \cH^{\max, d},
    \end{equation}
    where $R^d \PP^3$ is the smoothable component of the Hilbert scheme of
    points of $\PP^3$.
    \begin{lemma}\label{ref:tangentlemma:lem}
        The map $\varphi_d$ is dominating for all $8\leq d\leq 95$.
    \end{lemma}
    \begin{proof}
        First, we prove that for every $8\leq d\leq 95$ there is a smooth point $x\in R^{d} \PP^3$ such that the tangent map
        \begin{equation*}
            T_{\varphi_{d}}: \left(T R^{d} \PP^3\right) _{x}
            \to \left(T\cH^{\max, d}\right) _{\varphi_{d}(x)}
        \end{equation*}
        is surjective.
        This is verified by a direct computer calculation.
        See the accompanying package
        \emph{CombalggeomApprenticeshipsHilbert.m2}.
        Then by \cite[Theorem~17.11.1d, p.~83]{ega4-4} the morphism $f$ is smooth at
        $x$, thus flat, thus open, and thus the claim.
        \qed
    \end{proof}

    \begin{proposition}\label{ref:1361final:prop}
        For all $d\leq 95$ all schemes in $\cH^{\max, d}$ are
        smoothable.
    \end{proposition}
    \begin{proof}
        By Lemma~\ref{ref:upto95irreducibility}, the locus $\cH^{\max, d}$ is irreducible.
        For $d < 8$ all schemes are smoothable by~\cite{CEVV}. Assume $d\geq
        8$.
        By Lemma~\ref{ref:tangentlemma:lem}, the map~$\varphi_{d}$
        from~\eqref{eq:dominating} is dominating. Hence, a general element of
        $\cH^{\max, d}$ is smoothable. But smoothability is a closed
        property, thus, all elements of $\cH^{\max, d}$ are smoothable.
    \qed
    \end{proof}

    \begin{proof}[Proof of Proposition~\ref{ref:upto95:prop}]
        When $d\leq 95$ the claim follows from
        Proposition~\ref{ref:1361final:prop}. When $d\geq 96$ the claim
        follows by dimension count, as in~\cite{iarrobino_reducibility}.
    \qed
    \end{proof}

    \begin{remark}[Comparison with the case of $8$ points in $A^4$]\label{ref:furtherwork:remark}
        From the case $d = 96$ onwards we do not get
        a surjective tangent map and, indeed, the dimension of the family
        $\cH^{\max, 96} = \cH_3^{(1, 3, 6, 10, 20, 35, 21)}$ is equal to $3\cdot 96$,
        thus, a general member of this family cannot be smoothable for
        dimensional reasons.
        (Points in $\cH^{\max,96}$ define schemes supported at a single point,
        so $\cH^{\max,96}$ would have to be contained in the boundary of $R^{96}_3$.)
        This was the original example
        of~\cite{iarrobino_reducibility}.
        Our methods show that $\cH^{\max, d}$ is smoothable for all $d
        \leq 95$, hence, the bound $d = 96$ obtained
        in~\cite{iarrobino_reducibility} is sharp for this method. Note that
        in~\cite{iarrobino_compressed_artin} another, only partially related,
        method was used to prove that $H^{d}_3$ is reducible for $d\geq 78$.
        It is currently unclear whether this other method can yield
        irreducible components for $d\leq 77$.

        Even though $T_{\varphi_d}$ is not surjective for $d\geq 96$, we conjecture that the
        maps $T_{\varphi_{d}}$ are of maximal rank.
        This is no longer true for $\A^4$: in fact $T_{\varphi_{8}\,
        \A^4}$ has $20$-dimensional image in the $21$-dimensional
        Grassmannian $\Gr(3,10)$, which accounts for the fact that there are
        nonsmoothable ideals of degree $8$ in $\A^4$, as proven in~\cite{CEVV}.
        An explicit example of such a scheme in
        $\A^4 = \Spec \kk[\alpha, \beta,\gamma, \delta]$ is given by
        the ideal
        $(\alpha^2, \alpha\beta, \beta^2, \alpha\delta + \beta\gamma,
        \gamma^2, \gamma\delta, \delta^2) = \langle xz,xw,yz,yw,xy-zw \rangle^\perp$, see~\cite[Proposition~5.1]{CEVV}.
        This scheme gives an answer to
        \cite[Problem 3 on Parameters and Moduli]{Bernd}.
    \end{remark}


\begin{acknowledgement}
This article was initiated during the Apprenticeship
     Weeks (22 August--2 September 2016), led by Bernd
     Sturmfels, as part of the Combinatorial Algebraic
     Geometry Semester at the Fields Institute.
The authors wish to thank the Fields Institute,
the organizers of the Thematic Program on Combinatorial Algebraic Geometry in Fall 2016,
and the organizers of the Apprenticeship Weeks which took place during the program.
We are very grateful to Mark Huibregtse, Anthony Iarrobino, Gary Kennedy,
Greg Smith, Bernd Sturmfels, and several anonymous referees
for numerous helpful
comments.
This work was supported by a grant from the Simons Foundation (\#354574, Zach Teitler).
JJ was supported by Polish National Science Center, project
2014/13/N/ST1/02640. BIUN was supported by NRC project 144013.
\end{acknowledgement}

\def\cdprime{$''$}

\end{document}